\documentclass[12pt]{article}
\usepackage[centertags]{amsmath}
\usepackage{amsfonts}
\usepackage{amssymb}
\usepackage{amsthm}
\usepackage[top=3cm, bottom=3cm,left=3cm, right=3cm]{geometry}
\begin{document}
 \newtheorem{theorem}{Theorem}[section]
 \newtheorem{corollary}[theorem]{Corollary}
 \newtheorem{lemma}[theorem]{Lemma}
 \newtheorem{proposition}[theorem]{Proposition}
 \theoremstyle{definition}
 \newtheorem{definition}[theorem]{Definition}
 \theoremstyle{remark}
 \newtheorem{remark}[theorem]{Remark}
\newtheorem{example}[theorem]{Example}
 \numberwithin{equation}{section}

\makeatletter
\newcommand{\RM}{\mathbb{R}}
\newcommand{\ZM}{\mathbb{Z}}
\newcommand{\QM}{\mathbb{Q}}
\newcommand{\NM}{\mathbb{N}}
\newcommand{\CM}{\mathbb{C}}
\newcommand{\Mt}{\mathop{\rm Mult}\nolimits}

\title{On Some Lie Bialgebra Structures on Polynomial Algebras and their
Quantization}
\author{S. M. Khoroshkin$^{1}$, I. I. Pop$^{2}$, M. E. Samsonov$^{3}$,
A. A. Stolin$^{4}$, V. N. Tolstoy$^{5}$} 
\date{}
\maketitle
\begin{center}
\footnotesize{$^{1}$Institute for Theoretical and 
Experimental Physics, Moscow, Russia.\\
$^{2}$Department of Mathematical Sciences, G\"oteborg University, 
G\"oteborg, Sweden.\\
$^{3}$Dublin Institute for Advanced Studies, Dublin, Ireland.\\
$^{4}$Department of Mathematical Sciences, G\"oteborg University, 
G\"oteborg, Sweden.\\
$^{5}$Institute for Nuclear Physics, Moscow State University, Moscow, Russia.}
\end{center}

\begin{abstract}
We study classical twists of Lie bialgebra structures on the polynomial current algebra $\mathfrak{g}[u]$, where
$\mathfrak{g}$ is a simple complex finite-dimensional Lie algebra. 
We focus on the structures induced by the so-called quasi-trigonometric solutions of the classical
Yang-Baxter equation. 
It turns out that quasi-trigonometric $r$-matrices fall into classes labelled by the vertices of the extended Dynkin diagram of $\mathfrak{g}$. We give complete classification of quasi-trigonometric $r$-matrices belonging
to multiplicity free simple roots (which have coefficient 1 in the decomposition of the maximal root).
We quantize solutions corresponding to the first root of $\mathfrak{sl}(n)$. 
\end{abstract}

\section{Introduction}
Recall that, given a Lie algebra $\mathfrak{g}$, the classical Yang-Baxter equation
(CYBE) with one spectral parameter is the equation
\begin{equation}
[X^{12}(u),X^{13}(u+v)]+[X^{12}(u),X^{23}(v)]+[X^{13}(u+v),X^{23}(v)]=0,\label{eq:CYBE}
\end{equation}
where $X(u)$ is a meromorphic function of one complex variable $u$, defined in a
neighbourhood of $0$, taking values in $\mathfrak{g}\otimes\mathfrak{g}$. In their
outstanding paper \cite{BD2}, A. Belavin and V. Drinfeld investigated solutions of the
CYBE for a simple complex Lie algebra $\mathfrak{g}$. They considered so-called
\emph{nondegenerate solutions} (i.e. $X(u)$ has maximal rank for generic $u$). It was
proved in \cite{BD2} that nondegenerate solutions are of three types: rational,
trigonometric and elliptic. Moreover the authors completely classified trigonometric and
elliptic solutions, the last ones for the case $\mathbf{g}=\mathfrak{sl}(n)$.

One can see that any rational solution of CYBE provides the Lie bialgebra structure on
polynomial Lie algebra $\mathfrak{g}[u]$ for a simple Lie algebra $\mathfrak{g}$. On the
contrary, there are no clear Lie bialgebra structures related to elliptic solutions of
CYBE.

For trigonometric solutions of CYBE, the situation is as follows. Any trigonometric
solution has the form $Y(e^{k(u-v)})$, where $Y$ is a $\mathfrak{g}\otimes
\mathfrak{g}$-valued rational function and $k$ is some constant. After setting $e^{ku}=z$
this solution provides a Lie bialgebra structure either on  Lie algebra $\mathfrak{g}[z,z
^{-1}]$ or on its twisted version but does not induce, generally speaking, a Lie
bialgebra structure on the polynomial Lie algebra $\mathfrak{g}[z]$.

Therefore we are motivated to introduce a  class of solutions of 'trigonometric' type
that will induce Lie bialgebra structures on $\mathfrak{g}[u]$. Let $\Omega$ denote the
quadratic Casimir element of $\mathfrak{g}$. We say that a solution $X$ of the CYBE is
\emph{quasi-trigonometric} if it is of the form:
\begin{equation}
X(u,v)=\frac{v\Omega}{u-v}+p(u,v),\label{eq:trig}
\end{equation}
where $p(u,v)$ is a polynomial with coefficients in $\mathfrak{g}\otimes\mathfrak{g}$. We
will prove that by applying a certain holomorphic transformation and a change of
variables, any quasi-trigonometric solution becomes trigonometric, in the sense of
Belavin-Drinfeld classification.

In the works \cite{D1,D2} to any Lie bialgebra V. Drinfeld assigned another Lie bialgebra,
the so-called classical double. F. Montaner and E. Zelmanov \cite{MZ} proved that for any
Lie bialgebra structure on $\mathfrak{g}[u]$ its classical double is isomorphic as a Lie
algebra to one of four Lie algebras. We will consider two of them:
$\mathfrak{g}((u^{-1}))$ and $\mathfrak{g}((u^{-1}))\oplus\mathfrak{g}$.

The study of the Lie bialgebra structures given by quasi-trigonometric solutions will be
based on the description of the classical double. We show that all quasi-trigonometric
solutions induce the same classical double $\mathfrak{g}((u^{-1}))\oplus\mathfrak{g}$.
Moreover we construct a one-to-one correspondence between this type of solutions and a
special class of Lagrangian subalgebras of the $\mathfrak{g}((u^{-1}))\oplus
\mathfrak{g}$. It turns out that such Lagrangian subalgebras can be embedded into some
maximal orders of $\mathfrak{g}((u^{-1}))\oplus\mathfrak{g}$, which correspond to
vertices of the extended Dynkin diagram of $\mathfrak{g}$. This embedding enables us to
classify quasi-trigonometric solutions of CYBE which correspond to multiplicity free roots. 
We also use the classification of Manin triples for
reductive Lie algebras in terms of generalized Belavin-Drinfeld data obtained by  P.
Delorme \cite{Del}. In particular, we get a complete combinatorial description of
quasi-trigonometric solutions of CYBE, related to Lie algebra $\mathfrak{sl}(n)$.

The goal of the second part of the paper is to propose a quantization scheme for some of the
Lie bialgebra structures on $\mathfrak{g}[u]$ for $\mathfrak{g}=\mathfrak{sl}(n)$
described in the first part of the paper.  In all these cases the quantization is given
by an explicit construction of the corresponding twist. More precisely, the corresponding
Hopf algebra is isomorphic to $U_q(\mathfrak{g}[u])$ with twisted comultiplication, where
$U_q(\mathfrak{g}[u])$ is defined as certain subalgebra of quantum affine algebra $U_q(
\widehat{\mathfrak{g}})$. This result confirms the natural conjecture made
in \cite{KPST} and recently proved in \cite{Hal}: any classical twist can be quantized.

For the construction of twist quantizations of quasi-trigonometric solutions of CYBE, we
use nontrivial embedings of certain Hopf subalgebras of the quantized universal
enveloping algebra $U_q( {\mathfrak{sl}}_{n+1})$, called seaweed algebras \cite{DK} into
$U_q( \widehat{\mathfrak{sl}_n})$. This enables us to 'affinize' the finite - dimensional twists constructed in \cite{ESS} and \cite{IO}.

\section{Lie bialgebra structures and classical twists}

Let $\mathfrak{g}$ denote an arbitrary complex Lie algebra. We recall that a Lie
bialgebra structure on $\mathfrak{g}$ is a 1-cocycle $\delta:\mathfrak{g}
\longrightarrow\wedge^{2}\mathfrak{g}$ which satisfies the co-Jacobi identity. In other
words, $\delta$ provides a Lie algebra structure for $\mathfrak{g}^{*}$ compatible with
the structure of $\mathfrak{g}$.

To any Lie bialgebra $(\mathfrak{g},\delta)$ one associates the so-called \emph{classical
double} $D(\mathfrak{g},\delta)$. It is defined as the unique Lie algebra structure on
the vector space $\mathfrak{g}\oplus\mathfrak{g}^{*}$ such that:

a) it induces the given Lie algebra structures on $\mathfrak{g}$ and $\mathfrak{g}^{*}$

b) the bilinear form $Q$ defined by
\begin{equation}
Q(x_{1}+l_{1},x_{2}+l_{2})=l_{1}(x_{2})+l_{2}(x_{1})
\end{equation}
 is invariant with respect to the adjoint representation of
$\mathfrak{g}\oplus\mathfrak{g}^{*}$.

Let $\delta_{1}$ be a Lie bialgebra structure on $\mathfrak{g}$. Suppose
$s\in\wedge^{2}\mathfrak{g}$ satisfies
\begin{equation}
[s^{12},s^{13}]+[s^{12},s^{23}]+[s^{13},s^{23}]=\rm{Alt}(\delta_{1}\otimes \rm{id})(s),
\end{equation}
where $\mathrm{Alt}(x):=x^{123}+x^{231}+x^{312}$ for any $x\in\mathfrak{g}^{\otimes3}$. Then
\begin{equation}
\delta_{2}(a):=\delta_{1}(a)+[a\otimes1+1\otimes a,s]
\end{equation}
defines a Lie bialgebra structure on $\mathfrak{g}$. We call $s$ a \emph{classical twist}
and say that the bialgebra structures $(\mathfrak{g},\delta_{1})$ and $(\mathfrak{g},
\delta_{2})$ are \emph{related by a classical twist}.

The construction of the double suggests another notion of equivalence between bialgebras.
Namely, we say that Lie bialgebra structures $\delta_{1}$ and $\delta_{2}$ on
$\mathfrak{g}$ are   \emph{in the same twisting class} if there is a Lie algebra
isomorphism $f:D(\mathfrak{g},\delta_{1})\longrightarrow D(\mathfrak{g},\delta_{2})$
satisfying the properties:

1) $Q_{1}(x,y)=Q_{2}(f(x),f(y))$ for any $x,y \in D(\mathfrak{g},\delta_{1})$, where
$Q_{i}$ denotes the canonical form on $D(\mathfrak{g},\delta_{i})$, $i=1,2$.

2) $f\circ j_{1}=j_{2}$, where $j_{i}$ is the canonical embedding of $\mathfrak{g}$ in
$D(\mathfrak{g},\delta_{i})$.

For a finite-dimensional $\mathfrak{g}$, it was shown in \cite{KaS} that two Lie
bialgebra structures are in the same twisting class if and only if they are
related by a classical twist.

\begin{example}
\label{ex1}Let $\mathfrak{g}$ be finite-dimensional. All Lie bialgebra structures induced
by triangular $r$-matrices are related by classical twists. The classical double
corresponding to any triangular $r$-matrix is isomorphic to the semidirect sum
$\mathfrak{g}\dotplus\mathfrak{g}^{*}$ such that $\mathfrak{g}^{*}$ is a commutative
ideal and $[a,l]=\mathrm{ad} ^{*}(a)(l)$ for any $a\in\mathfrak{g}$ and
$l\in\mathfrak{g}^{*}$.
\end{example}
Another example of twisting is the following:

\begin{example}
\label{ex2}Suppose $\mathfrak{g}$ is simple and let $\delta_{0}$ be the Lie bialgebra
structure induced by the standard Drinfeld-Jimbo $r$-matrix. Then the entire
Belavin-Drinfeld list \cite{BD1} is obtained by twisting the standard structure
$\delta_{0}$. The classical double corresponding to any $r$-matrix from this list is
isomorphic to $\mathfrak{g}\oplus\mathfrak{g}$.
\end{example}
Now, if we pass to the case of infinite-dimensional Lie bialgebra structures, we
encounter more examples of twisting.

Let us recall several facts from the theory of rational solutions as it was developed in
\cite{S1}. We let again $\mathfrak{g}$ denote a simple Lie algebra. Denote by $K$ the
Killing form and let $\Omega$ be the corresponding Casimir element of $\mathfrak{g}$. We
look for functions $X:\mathbb{C}^{2}\longrightarrow\mathfrak{g}\otimes\mathfrak{g}$ which
satisfy
\begin{equation}
[X^{12}(u_{1},u_{2}),X^{13}(u_{1},u_{3})]+[X^{12}(u_{1},u_{2}),X^{23}(u_{2},u_{3})]+
\label{eq:4}
\end{equation}
\[
+[X^{13}(u_{1},u_{3}),X^{23}(u_{2},u_{3})]=0,\]
\begin{equation}
X^{12}(u,v)=-X^{21}(v,u).\label{eq:5}
\end{equation}

\begin{remark}
We will call these two equations \emph{the classical Yang-Baxter equation (CYBE)}. In the
case of rational and quasi-trigonometric solutions, the unitarity condition (\ref{eq:5})
can actually be dropped. We will prove in Appendix that (\ref{eq:5}) is a consequence of
(\ref{eq:4}).
\end{remark}

\begin{definition}
Let $X(u,v)=\frac{\Omega}{u-v}+p(u,v)$ be a function from $\mathbb{C}^{2}$ to
$\mathfrak{g}\otimes\mathfrak{g}$, where $p(u,v)$ is a polynomial with coefficients in
$\mathfrak{g}\otimes\mathfrak{g}$. If $X$ satisfies the CYBE, we say that $X$ is a
\emph{rational} solution.

Two rational solutions $X_{1}$ and $X_{2}$ are called \emph{gauge equivalent} if there
exists $\sigma(u)\in \mathrm{Aut}(\mathfrak{g}[u])$ 
such that $X_{2}(u,v)=(\sigma(u)\otimes\sigma(v))X_{1}(u,v)$, 
where $\mathrm{Aut}(\mathfrak{g}[u])$ denotes the group of automorphisms of
$\mathfrak{g}[u]$ considered as an algebra over $\mathbb{C}[u]$.
\end{definition}
\begin{remark}
\label{rem:D} It was proved in \cite{S1} that any rational solution can be brought by
means of a gauge transformation to the form:
\[
X(u,v)=\frac{\Omega}{u-v}+p_{00}+p_{10}u+p_{01}v+p_{11}uv,
\]
where $p_{00}$, $p_{10}$, $p_{01}$, $p_{11}$$\in\mathfrak{g}\otimes\mathfrak{g}$.
\end{remark}
We recall that any rational solution induces a Lie bialgebra structure on the polynomial
current algebra $\mathfrak{g}[u]$. Let us consider a rational solution $X$ and define the
map $\delta_{X}:\mathfrak{g}[u]\longrightarrow\mathfrak{g}[u]\wedge\mathfrak{g}[u]$ by
\begin{equation}
\delta_{X}(a(u))=[X(u,v),a(u)\otimes1+1\otimes a(v)],\label{eq:7}
\end{equation}
for any $a(u)\in\mathfrak{g}[u]$. Obviously $\delta_{X}$ is a 1-cocycle and therefore
induces a Lie bialgebra structure on $\mathfrak{g}[u]$.

The following result, proved in \cite{S1}, shows that all Lie bialgebra structures
corresponding to rational solutions have the same classical double.

Let $\mathbb{C}[[u^{-1}]]$ be the ring of formal power series in $u^{-1}$ and
$\mathbb{C}((u^{-1}))$ its field of quotients. Consider the Lie algebras
$\mathfrak{g}[u]=\mathfrak{g}\otimes\mathbb{C}[u]$,
$\mathfrak{g}[[u^{-1}]]=\mathfrak{g}\otimes\mathbb{C}[[u^{-1}]]$ and
$\mathfrak{g}(u^{-1}))=\mathfrak{g}\otimes\mathbb{C}((u^{-1}))$.

Let $D_{X}(\mathfrak{g}[u])$ be the classical double corresponding to a rational solution
$X$ of the CYBE. Then $D_{X}(\mathfrak{g}[u])$ and $\mathfrak{g}((u^{-1}))$ are
isomorphic as Lie algebras, with inner product which has the following form on
$\mathfrak{g}((u^{-1}))$:
\begin{equation}
Q(f(u),g(u))=Res_{u=0}K(f(u),g(u)),\label{eq:8}
\end{equation}
where $f(u),g(u)\in\mathfrak{g}((u^{-1}))$.

\begin{remark}
This result means in fact that $\mathfrak{g}((u^{-1}))$ can be represented as a Manin triple
$\mathfrak{g}((u^{-1}))=\mathfrak{g}[u]\oplus W$, where $W$ is a Lagrangian subalgebra
with respect to the invariant form $Q$.
\end{remark}
In the case $\mathfrak{g}=\mathfrak{sl}(n)$ all rational solutions were described in
the following way:

Let $d_{k}=diag(1,...,1,u,...,u)$ ($k$ many 1's), $0\leq k\leq[\frac{n}{2}]$. Then it was
proved in \cite{S1} that every rational solution of the CYBE defines some Lagrangian
subalgebra $W$ contained in $d_{k}^{-1}\mathfrak{sl}(n)[[u^{-1}]]d_{k}$ for some $k$.
These subalgebras are in one-to-one correspondence with pairs $(L,B)$ verifying:

(1) $L$ is a subalgebra of $\mathfrak{sl}(n)$ such that $L+P_{k}=\mathfrak{sl}(n)$, where
$P_{k}$ denotes the maximal parabolic subalgebra of $\mathfrak{sl}(n)$ not containing the
root vector $e_{\alpha_{k}}$of the simple root $\alpha_{k}$;

(2) $B$ is a 2-cocycle on $L$ which is nondegenerate on $L\cap P_{k}$.

In case $\mathfrak{g}=\mathfrak{sl}(2)$ one has just two non-standard rational
$r$-matrices, up to gauge equivalence:
\begin{equation}
X_{1}(u,v)=\frac{\Omega}{u-v}+h_{\alpha}\wedge e_{-\alpha}\label{eq:R1}
\end{equation}
and
\begin{equation}
X_{2}(u,v)=\frac{\Omega}{u-v}+ue_{-\alpha}\otimes h_{\alpha}-vh_{\alpha}\otimes
e_{-\alpha},\label{eq:R2}
\end{equation}
where $e_{\alpha}=e_{12}$, $e_{-\alpha}=e_{21}$ and $h_{\alpha}=e_{11}-e_{22}$ is the
usual basis of $\mathfrak{sl}(2)$.

\section{Quasi-trigonometric solutions of the CYBE}
Another interesting case of infinite-dimensional Lie bialgebra structures on
$\mathfrak{g}[u]$ is provided by a class of trigonometric type solutions of the CYBE,
called \emph{quasi-trigonometric solutions}. These solutions were first introduced in
\cite{KPST}.
\begin{definition}
We say that a solution $X$ of the CYBE is \emph{quasi-trigonometric} if it is of the
form: \begin{equation} X(u,v)=\frac{v\Omega}{u-v}+p(u,v),\label{eq:10b}\end{equation}
 where $p(u,v)$ is a polynomial with coefficients in $\mathfrak{g}\otimes\mathfrak{g}$.
\end{definition}
The term \emph{quasi-trigonometric} is motivated by the relationship between this  type
of solutions of CYBE and trigonometric solutions in the Belavin-Drinfeld meaning. The
following result, whose proof we give in the Appendix, illustrates this fact.
\begin{theorem}
Let $X(u,v)$ be a quasi-trigonometric solution of the CYBE. There exists a holomorphic
transformation and a change of variables such that $X(u,v)$ becomes a trigonometric
solution, in the sense of Belavin-Drinfeld classification.
\end{theorem}
\begin{example}
A function $X(u,v)=\frac{v\Omega}{u-v}+r$, where $r\in\mathfrak{g}\otimes\mathfrak{g}$,
satisfies the CYBE if and only if $r$ is a solution of the \emph{modified classical
Yang-Baxter equation}, i.e.
\begin{equation}
r+r^{21}=\Omega
\end{equation}
\begin{equation}
[r^{12},r^{13}]+[r^{12},r^{23}]+[r^{13},r^{23}]=0.
\end{equation}
Let $r_{0}$ denote the standard Drinfeld-Jimbo r-matrix. We fix a Cartan subalgebra
$\mathfrak{h}$ and the associated root system. We choose a system of generators
$e_{\alpha}$, $e_{-\alpha}$ and $h_{\alpha}$, where $\alpha$ is a positive root, such
that $K(e_{\alpha},e_{-\alpha})=1$. Then
\begin{equation}
r_{0}=\frac{1}{2}(\sum_{\alpha>0}e_{\alpha}\wedge e_{-\alpha}+\Omega), \label{eq:DJ}
\end{equation}
Correspondingly we have a quasi-trigonometric solution
\begin{equation}\label{X_{0}}
X_{0}(u,v)=\frac{v\Omega}{u-v}+r_{0}.
\end{equation}
 \end{example}
\begin{definition}
A quasi-trigonometric solution $X(u,v)=\frac{v\Omega}{u-v}+p(u,v)$ is called
\emph{quasi-constant} if $p(u,v)$ is a constant polynomial.
\end{definition}
\begin{proposition}
Let $\mathrm{Aut}(\mathfrak{g}[u])$ denote the group of automorphisms of $\mathfrak{g}[u]$
considered as an algebra over $\mathbb{C}[u]$. Let $X_{1}$ be a quasi-trigonometric
solution and $\sigma(u)\in \mathrm{Aut}(\mathfrak{g}[u])$. Then
\[
X_{2}(u,v)=(\sigma(u)\otimes\sigma(v))X_{1}(u,v)
\]
is also a quasi-trigonometric solution.
\end{proposition}
\begin{proof}
Let $X_{1}(u,v)=\frac{v\Omega}{u-v}+p(u,v)$. Since $X_{2}$ obviously satisfies the CYBE,
it is enough to check that $X_{2}$ is quasi-trigonometric. We have the following: \[
X_{2}(u,v)=\left(\frac{v(\sigma(u)-\sigma(v))}{u-v}\otimes\sigma(v)\right)
\Omega+\frac{v}{u-v}(\sigma(v)\otimes\sigma(v))\Omega+\]
\[
+(\sigma(u)\otimes\sigma(v))p(u,v).\] Let
$p_{1}(u,v):=(\frac{v(\sigma(u)-\sigma(v))}{u-v}\otimes\sigma(v))\Omega$ and
$p_{2}(u,v):=(\sigma(u)\otimes\sigma(v))p(u,v)$. These are polynomial functions in $u$,
$v$. Since $(\sigma(v)\otimes\sigma(v))\Omega=\Omega$, we obtain
\[
X_{2}(u,v)=\frac{v\Omega}{u-v}+p_{1}(u,v)+p_{2}(u,v)
\]
and this ends the proof.
\end{proof}
\begin{definition}
Two quasi-trigonometric solutions $X_{1}$ and $X_{2}$ are called \emph{gauge equivalent}
if there exists $\sigma(u)\in \mathrm{Aut}(\mathfrak{g}[u])$ such that
\begin{equation}
X_{2}(u,v)=(\sigma(u)\otimes\sigma(v))X_{1}(u,v).\label{G2}
\end{equation}

\end{definition}
Any quasi-trigonometric solution $X$ of the CYBE induces a Lie bialgebra structure on
$\mathfrak{g}[u]$. Let $\delta_{X}$ be the 1-cocycle defined by
\begin{equation}
\delta_{X}(a(u))=[X(u,v),a(u)\otimes1+1\otimes a(v)],\label{eq:11}
\end{equation}
for any $a(u)\in\mathfrak{g}[u]$.

It is expected that all Lie bialgebra structures corresponding to quasi-trigonometric
solutions induce the same classical double.

Let us consider the direct sum of Lie algebras
$\mathfrak{g}((u^{-1}))\oplus\mathfrak{g}$, together with the following invariant
bilinear form:
\begin{equation}
Q((f(u),a),(g(u),b))=K(f(u),g(u))_{0}-K(a,b).\label{eq:12}
\end{equation}
Here the index zero means that we have taken the free term in the series expansion.
\begin{remark}\label{rem:id}
The Lie algebra $\mathfrak{g}[u]$ is embedded into $\mathfrak{g}((u^{-1}))
\oplus\mathfrak{g}$ via $a(u)\longmapsto(a(u),a(0))$ and is naturally identified with
\begin{equation}
V_{0}:=\{(a(u),a(0));a(u)\in\mathfrak{g}[u]\}.\label{eq:V_0}
\end{equation}

Consider the following Lie subalgebra of $\mathfrak{g}((u^{-1}))\oplus\mathfrak{g}$:
\begin{equation}\label{W0}
W_{0}=\{(a+f(z),b):f\in
z^{-1}\mathfrak{g}[[z^{-1}]],a\in\mathfrak{b}_{+},b\in\mathfrak{b}_{-},a_{\mathfrak{h}}+b_{\mathfrak{h}}=0\}.
\end{equation}

Here $\mathfrak{h}$ is the fixed Cartan subalgebra of $\mathfrak{g}$,
$\mathfrak{b_{\pm}}$ are the positive (negative) Borel subalgebras and $a_{\mathfrak{h}}$
denotes the Cartan part of $a$. We make the remark that $V_{0}\oplus
W_{0}=\mathfrak{g}((u^{-1}))\oplus\mathfrak{g}$ and both $V_{0}$ and $W_{0}$ are
isotropic with respect to the form $Q$.
\end{remark}

In order to show that all quasi-trigonometric solutions induce the same classical double,
we will first prove the following result:

\begin{theorem}
\label{thm:4}There exists a natural one-to-one correspondence between quasi - trigonometric
solutions of the CYBE and linear subspaces $W$ of
$\mathfrak{g}((u^{-1}))\oplus\mathfrak{g}$ such that

1) $W$ is a Lie subalgebra in $\mathfrak{g}((u^{-1}))\oplus\mathfrak{g}$ such that
$W\supseteq u^{-N}\mathfrak{g}[[u^{-1}]]$ for some $N>0$;

2) $W\oplus V_{0}=\mathfrak{g}((u^{-1}))\oplus\mathfrak{g}$;

3) $W$ is a Lagrangian subspace with respect to the inner product of $\mathfrak{g}
((u^{-1}))\oplus\mathfrak{g}$.
\end{theorem}
\begin{proof}
Let $V_{0}$ and $W_{0}$ be the Lie algebras given in Remark \ref{rem:id}. We choose dual
bases in $V_{0}$ and $W_{0}$ respectively. Let $\{ k_{j}\}$ be an orthonormal basis in
$\mathfrak{h}$. The canonical basis of $V_{0}$ is formed by $e_{\alpha}u^{k}$,
$e_{-\alpha}u^{k}$, $k_{j}u^{k}$ for any $\alpha>0$, $k>0$ and all $j$;
$(e_{-\alpha},e_{-\alpha})$, $(e_{\alpha},e_{\alpha})$ for any $\alpha>0$, and
$(k_{j},k_{j})$, for all $j$. The dual basis of $W_{0}$ is the following:
$e_{-\alpha}u^{-k}$, $e_{\alpha}u^{-k}$, $k_{j}u^{-k}$ for any $\alpha>0$, $k>0$ and all
$j$; $(e_{\alpha},0)$ and $(0,-e_{-\alpha})$ for all $\alpha>0$, and
$\frac{1}{2}(k_{j},-k_{j})$, for all $j$. Let us simply denote these dual bases by $\{
v_{i}\}$ and $\{ w_{0}^{i}\}$ respectively. We notice that the quasi-trigonometric
solution $X_{0}$ can be written as
\begin{equation}
X_{0}(u,v)=(\tau\otimes\tau)(\sum_{i}w_{0}^{i}\otimes v_{i}),\label{eq:14'}
\end{equation}
where $\tau$ denotes the projection of $\mathfrak{g}((u^{-1}))\oplus\mathfrak{g}$ onto
$\mathfrak{g}((u^{-1}))$.

We denote by $Hom_{c}(W_{0},V_{0})$ the space of those linear maps
$F:W_{0}\longrightarrow V_{0}$ such that $KerF\supseteq u^{-N}\mathfrak{g}[[u^{-1}]]$ for
some $N>0$. It is the space of linear maps $F$ which are continuous with respect to the
``$u^{-1}$ - adic'' topology.

Let us contruct a linear isomorphism $\Phi:V_{0}\otimes V_{0}\longrightarrow
Hom_{c}(W_{0},V_{0})$ in the following way: \begin{equation} \Phi(x\otimes
y)(w_{0})=Q(w_{0},y)\cdot x,\label{eq:fi}\end{equation} for any $x$, $y$ $\in V_{0}$ and
any $w_{0}\in W_{0}$. It is easy to check that $\Phi$ is indeed an isomorphism. The
inverse mapping is $\Psi:Hom_{cont}(W_{0},V_{0})\longrightarrow V_{0}\otimes V_{0}$
defined by
\begin{equation}
\Psi(F)=\sum_{i}F(w_{0}^{i})\otimes v_{i}.\label{eq:psi}
\end{equation}
We make the remark that this sum is finite since $F(w_{0}^{i})\neq0$ only for a finite
number of indices $i$.

The next step is to construct a bijection between $Hom_{c}(W_{0},V_{0})$ and the set $L$
of linear subspaces $W$ of $\mathfrak{g}((u^{-1}))\oplus\mathfrak{g}$ such that $W\oplus
V_{0}=\mathfrak{g}((u^{-1}))\oplus\mathfrak{g}$ and $W\supseteq u^{-N}\mathbf{g}
[[u^{-1}]]$ for some $N>0$. This can be done in a very natural way. For any $F\in
Hom_{cont}(W_{0},V_{0})$ we take \begin{equation} W(F)=\{ w_{0}+F(w_{0});w_{0}\in
W_{0}\}.\label{eq:14b}\end{equation} The inverse mapping associates to any $W$ the linear
function $F_{W}$ such that for any $w_{0}\in W_{0}$, $F_{W}(w_{0})=-v$, uniquely defined
by the decomposition $w_{0}=w+v_{0}$ with $w\in W$ and $v_{0}\in V_{0}$.

Therefore we have a bijection between $V_{0}\otimes V_{0}$ and $L$. By a straightforward
computation, one can show that a tensor $r(u,v)\in V_{0}\otimes V_{0}$ satisfies the
condition $r(u,v)=-r^{21}(v,u)$ if and only if the linear subspace $W(\Phi(r))$ is
Lagrangian with respect to $Q$.

Let us suppose now that $X(u,v)=X_{0}(u,v)+r(u,v)$ and $r(u,v)=-r^{21}(v,u)$. Then
$X(u,v)$ satisfies (\ref{eq:4}) if and only if $W(\Phi(r))$ is a Lie subalgebra of
$\mathbf{g}((u^{-1}))\oplus\mathbf{g}$. Indeed, since $r$ is unitary, we have that
$W(\Phi(r))$ is a Lagrangian subspace with respect to $Q$. It is enough to check that
$X(u,v)$ satisfies (\ref{eq:4}) if and only if
\begin{equation}
Q([w_{1}+\Phi(r)(w_{1}),w_{2}+\Phi(r)(w_{2})],w_{3}+\Phi(r)(w_{3}))=0
\label{eq:condition}
\end{equation}
for any elements $w_{1}$, $w_{2}$ and $w_{3}$ of $W_{0}$. This follows by direct
computations.

In conclusion, we see that a function $X(u,v)=\frac{v\Omega}{u-v}+p(u,v)$ is a
quasi-trigonometric solution if and only if $W(\Phi(p-r_{0}))$ is a Lagrangian subalgebra
of $\mathfrak{g}((u^{-1}))\oplus\mathfrak{g}$. This ends the proof.
\end{proof}
\begin{remark}
\label{rem:1}If $W$ is a Lagrangian subalgebra of $\mathfrak{g}((u^{-1}))\oplus
\mathfrak{g}$ satisfying the conditions of Theorem \ref{thm:4}, then the corresponding
solution $X(u,v)$ is constructed in the following way: take a basis $\{ w^{i}\}$ in $W$
which is dual to the canonical basis $\{ v_{i}\}$ of $V_{0}$ and construct the tensor
\begin{equation} \widetilde{r}(u,v)=\sum_{i}w^{i}\otimes
v_{i}.\label{eq:rtilda}\end{equation}

Let $\pi$ denote the projection of $\mathbf{g}((u^{-1}))\oplus\mathbf{g}$ onto
$\mathbf{g}[u]$ which is induced by the decomposition $\mathbf{g}((u^{-1}))\oplus
\mathbf{g}=V_{0}\oplus W_{0}$. Explicitly,
\begin{equation}
\pi(a_{n}z^{n}+...+a_{0}+a_{-1}u^{-1}+...,b)=a_{n}u^{n}+...+a_{1}u+ \label{eq:14'}
\end{equation}
\[
+\frac{1}{2}(a_{0\mathfrak{h}}+b_{\mathfrak{h}})+a_{0-}+b_{+}.\] Here
$a_{0}=a_{0\mathfrak{h}}+a_{0+}+a_{0-}$ and $b=b_{\mathfrak{h}}+b_{+}+b_{-}$ are the
decompositions with respect to $\mathbf{g}=\mathfrak{h}\oplus \mathfrak{n}_{+}\oplus
\mathfrak{n}_{-}$. Then
\begin{equation}
X(u,v)=X_{0}(u,v)+(\pi\otimes\pi)(\widetilde{r}(u,v)).\label{eq:14"}
\end{equation}
At this point we note the following fact that we will prove in Appendix:
\end{remark}
\begin{proposition}
Let $W$ be a Lie subalgebra satisfying conditions 2) and 3) of Theorem \ref{thm:4}. Let
$\widetilde{r}$ be constructed as in (\ref{eq:rtilda}). Assume $\widetilde{r}$ induces a
Lie bialgebra structure on $\mathbf{g}[u]$ by $\delta_{\widetilde{r}}(a(u))=
[\widetilde{r}(u,v),a(u)\otimes1+1\otimes a(v)]$. Then $W\supseteq u^{-N}\mathfrak{g}
[[u^{-1}]]$ for some positive $N$.
\end{proposition}
For any quasi-trigonometric solution $X$ of the CYBE denote by $D_{X}(\mathfrak{g}[u])$ the classical double of $\mathfrak{g}[u]$ corresponding to $X$.
\begin{theorem}\label{thdx} For any quasi-trigonometric solution $X$ of the CYBE there exists an isomorphism of Lie algebra $D_{X}(\mathfrak{g}[u])$
and $\mathfrak{g}((u^{-1}))\oplus\mathfrak{g}$ identical on $\mathfrak{g}[u]$, which
transforms the canonical bilinear form on $D_{X}(\mathfrak{g}[u])$ to the form $Q$.
\end{theorem}
\begin{proof}
One can easily check that if $W$ is a Lagrangian subalgebra of $\mathfrak{g}((u^{-1}))
\oplus\mathfrak{g}$, corresponding to a quasi-trigonometric solution
$X(u,v)=\frac{v\Omega}{u-v}+p(u,v)$, then $W$ is isomorphic to $(\mathfrak{g}[u])^{*}$
with the Lie algebra structure induced by $X$.

Indeed, with the notation introduced in Theorem \ref{thm:4}, let $F:=F_{W}$. It is enough
to check that for any $v_{0}\in V_{0}$ and $w_{1},w_{2}\in W_{0}$ the following equality
is satisfied:\[ Q(v_{0},[w_{1}+F(w_{1}),w_{2}+F(w_{2})])=<\delta_{X}(v_{0}),w_{1}\otimes
w_{2}>,\] where $<,>$ denotes the pairing between $V_{0}^{\otimes2}$ and
$W_{0}^{\otimes2}$ induced by $Q$. This equality is implied by the following identities:
\[Q(v_{0},[w_{1},w_{2}])=<\delta_{X_{0}}(v_{0}),w_{1}\otimes w_{2}>,\]
\[Q(v_{0},[F(w_{1}),w_{2}])=<[p-r_{0},1\otimes v_{0}],w_{1}\otimes w_{2}>,\]
\[Q(v_{0},[w_{1},F(w_{2})])=<[p-r_{0},v_{0}\otimes1],w_{1}\otimes
w_{2}>.\]
\end{proof}
\begin{remark} Theorem \ref{thdx} states in particular that all quasi-trigonometric
solutions of CYBE are in the same twisting class. This can be seen directly, since by the
definition of quasi-trigonometric solutions they are related to the solution
(\ref{X_{0}}) by classical twists.
\end{remark}

\begin{theorem}
\label{cor:gauge} Let $X_{1}$ and $X_{2}$ be quasi-trigonometric solutions of the CYBE.
Suppose that $W_{1}$ and $W_{2}$ are the corresponding Lagrangian subalgebras of
$\mathfrak{g}((u^{-1}))\oplus\mathfrak{g}$. Let $\sigma(u)\in \mathrm{Aut}(\mathfrak{g}[u]$) and
$\widetilde{\sigma}(u)$ be the automorphism of $\mathfrak{g}((u^{-1}))\oplus\mathfrak{g}$
induced by $\sigma(u)$. The following conditions are equivalent:

1) $X_{1}(u,v)=(\sigma(u)\otimes\sigma(v))X_{2}(u,v)$;

2) $W_{1}=\widetilde{\sigma}(u)W_{2}$.
\end{theorem}
\begin{proof}
1) $\Rightarrow$ 2). Let us begin by proving this for the particular case $X_{1}=X_{0}$
and $X_{2}=(\sigma(u)\otimes\sigma(v))X_{0}(u,v)$. The Lagrangian subalgebra
corresponding to $X_{0}$ is $W_{0}$ given by (\ref{W0}). On the other hand, one can
check the Lagrangian subalgebra $W_{2}$, corresponding to the solution $X_{2}$, consists
of elements
\[\widetilde{f}:=\sum_{i}(f,\widetilde{\sigma}(v_{i}))\cdot\widetilde{\sigma}(w_{0}^{i})=
\sum_{i}(\widetilde{\sigma}^{-1}(f),v_{i})\cdot\widetilde{\sigma}(w_{0}^{i}),\] for any
$f\in W_{0}$. Here $\{ v_{i}\}$ and $\{ w_{0}^{i}\}$ are the dual bases of $V_{0}$ and
$W_{0}$ introduced in the proof of Theorem \ref{thm:4}. We show that
$W_{2}=\widetilde{\sigma}(W_{0})$.

Let $g$ denote the projection of $\widetilde{\sigma}^{-1}(f)$ onto $W_{0}$ induced by the
decomposition $V_{0}\oplus W_{0}=\mathfrak{g}((u^{-1}))\oplus\mathfrak{g}$. Then \[g=
\sum_{i}(\widetilde{\sigma}^{-1}(f),v_{i})\cdot w_{0}^{i}\] which implies that
$\widetilde{f}=\widetilde{\sigma}(g)$. Therefore $W_{2}\subseteq
\widetilde{\sigma}(W_{0})$. The other inclusion is similar.

Let us pass to the general case. If $X_{1}(u,v)=X_{0}(u,v)+r(u,v)$ is a
quasi-trigonometric solution with $r(u,v)=\sum a_{k}u^{k}\otimes b_{j}v^{j}$, then the
corresponding $W_{1}$ consists of elements of the form $f+\sum(f,b_{j}u^{j})a_{k}u^{k}$,
for any $f$ in $W_{0}$. Now let $X_{2}(u,v)=(\sigma(u)\otimes\sigma(v))X_{1}(u,v)$. The
corresponding subalgebra $W_{2}$ is formed by elements of the form \[
\widetilde{f}_{r}:=\sum_{i}(f,\widetilde{\sigma}(v_{i}))\cdot\widetilde{\sigma}
(w_{0}^{i})+\sum(f,\widetilde{\sigma}(b_{j}u^{j}))\widetilde{\sigma}(a_{k}u^{k}).\] It is
easy to see that $\widetilde{f}_{r}=$$\widetilde{\sigma}(h)$, where
$h:=g+\sum(g,b_{j}u^{j})a_{k}u^{k}$ and $g$ is the projection of
$\widetilde{\sigma}^{-1}(f)$ onto $W_{0}$. These considerations prove that
$\widetilde{\sigma}(W_{1})=W_{2}$.

2) $\Rightarrow$ 1). Suppose that $W_{2}=\widetilde{\sigma}(W_{1})$. Let
$\widetilde{X}_{2}:=(\sigma(u)\otimes\sigma(v))X_{1}(u,v)$. It is a quasi-trigonometric
solution which has a corresponding Lagrangian subalgebra $\widetilde{W}_{2}$. Because 1)
$\Rightarrow$ 2) we obtain that $\widetilde{W}_{2}=\widetilde{\sigma}(W_{1})$ and thus
$W_{2}=\widetilde{W}_{2}$. Since the correspondence between solutions and subalgebras is
one-to-one, we get that $X_{2}=\widetilde{X}_{2}$.
\end{proof}
\begin{definition}
We will say that $W_{1}$ and $W_{2}$ are \emph{gauge equivalent} (with respect to
$\mathrm{Aut}(\mathfrak{g}[u])$) if condition 2) of Theorem \ref{cor:gauge} is satisfied.
\end{definition}
\begin{theorem}
\label{thm:5} Let $X(u,v)=\frac{v\Omega}{u-v}+p(u,v)$ be a quasi-trigonometric solution
of the CYBE and $W$ the corresponding Lagrangian subalgebra of $\mathfrak{g}((u^{-1}))
\oplus\mathfrak{g}$. Then the following are equivalent:

1) $p(u,v)$ is a constant polynomial.

2) $W$ is contained in $\mathfrak{g}[[u^{-1}]]\oplus\mathfrak{g}$.
\end{theorem}
\begin{proof}
We keep the notations from the proof of Theorem \ref{thm:4} and also those from Remark
\ref{rem:1}. Let $r(u,v)=p(u,v)-r_{0}$ and $F=\Phi(r(u,v))$. If $p(u,v)$ is constant,
then $F(w_{0})\in\mathfrak{g}\otimes \mathfrak{g}$ for any $w_{0}\in W_{0}$. Therefore
$W(F)\subseteq\mathfrak{g}[[u^{-1}]]\oplus \mathfrak{g}$. Conversely, let us suppose that
$W$ is included in $\mathfrak{g}[[u^{-1}]]\oplus\mathfrak{g}$. The orthogonal of
$\mathfrak{g}[[u^{-1}]]\oplus\mathfrak{g}$ with respect to $Q$ is obviously
$u^{-1}\mathfrak{g}[[u^{-1}]]$. Since $W$ is a Lagrangian subalgebra, it follows that $W$
contains $u^{-1}\mathfrak{g}[[u^{-1}]]$. According to the previous remark,
$r(u,v)=(\pi\otimes\pi)(\widetilde{r}(u,v))$, where $\pi$ is the projection onto
$\mathfrak{g}[u]$ induced by the decomposition
$\mathfrak{g}((u^{-1}))\oplus\mathfrak{g}=V_{0}\oplus W_{0}$ and
$\widetilde{r}(u,v)=\sum_{i}w^{i}\otimes v_{i}$. Now it is clear that $r$ is a constant.
This ends the proof.
\end{proof}

\section{Classification of quasi-trigonometric solutions}
We have seen that gauge equivalent solutions correspond to gauge equivalent subalgebras
$W$. Thus, the classification of quasi-trigonometric solutions is equivalent to the
classification of $W$ satisfying the conditions of Theorem \ref{thm:4}. In order to
classify such $W$ we will use a method from \cite{S3} which allows us to embed $W$ in a
suitable $\mathbb{C}$-subalgebra of $\mathfrak{g}((u^{-1}))\oplus\mathfrak{g}$.

Having fixed a Cartan subalgebra $\mathfrak{h}$ of $\mathfrak{g}$, let $R$ be the
correspondig set of roots and $\Gamma$ the set of simple roots. Denote by
$\mathfrak{g}_{\alpha}$ the root space corresponding to a root $\alpha$. Let
$\mathfrak{h}(\mathbb{R})$ be the set of all $h\in\mathfrak{h}$ such that
$\alpha(h)\in\mathbb{R}$ for all $\alpha\in R$. Consider the valuation on
$\mathbb{C}((u^{-1}))$ defined by $v(\sum_{k\geq n}a_{k}u^{-k})=n$. For any root $\alpha$
and any $h\in\mathfrak{h}(\mathbb{R})$, set $M_{\alpha}(h)$:=$\{f\in\mathbb{C}((u^{-1})):
v(f)\geq \alpha(h)\}$. Consider
\begin{equation}\label{eq3}
\mathbb{O}_{h}:=\mathfrak{h}[[u^{-1}]]\oplus(\oplus_{\alpha\in R}M_{\alpha}(h)
\otimes\mathfrak{g}_{\alpha}).
\end{equation}

As a direct corollary  of Theorem 4 from \cite{S3}, the following result can be deduced:
\begin{theorem}
Up to a gauge equivalence, any subalgebra $W$ which corresponds to a quasi-trigonometric
solution can be embedded into $\mathbb{O}_{h}\times\mathfrak{g}$, where $h$ is a vertex
of the following standard simplex $\{h\in\mathfrak{h}(\mathbb{R}):$ $\alpha(h)\geq 0$ for
all $\alpha\in\Gamma$ and $\alpha_{\max}\leq 1\}$.
\end{theorem}

Vertices of the above simplex correspond to vertices of the extended Dynkin diagram of
$\mathfrak{g}$, the correspondence being given by the following rule:
\[0\leftrightarrow\alpha_{\max}\]\[h_{i}\leftrightarrow\alpha_{i},\]
where $\alpha_{i}(h_{j})=\delta_{ij}/k_{j}$ and $k_{j}$ are given by the relation $\sum
k_{j}\alpha_{j}=\alpha_{\max}$. We will write $\mathbb{O}_{\alpha}$ instead of
$\mathbb{O}_{h}$ if $\alpha$ is the root which corresponds to the vertex $h$.
\begin{remark}
We have $\mathbb{O}_{\alpha_{\max}}=\mathfrak{g}[[u^{-1}]]$. We have already seen that a
quasi-trigonometric solution is quasi-constant if and only if its corresponding $W$ is embedded
into $\mathbb{O}_{\alpha_{\max}}\times\mathfrak{g}$.
\end{remark}
\begin{remark}
It might happen that two Lagrangian subalgebras $W_{1}$ and $W_{2}$ are gauge equivalent
even though they are embedded into different $\mathbb{O}_{\alpha_{1}}\times\mathfrak{g}$
and $\mathbb{O}_{\alpha_{2}}\times\mathfrak{g}$. If there exists an automorphism of the
Dynkin diagram of $\mathfrak{g}$ taking $\alpha_{1}$ into $\alpha_{2}$, then $W_{1}$ and
$W_{2}$ are gauge equivalent and the corresponding quasi-trigonometric solutions as well.
\end{remark}
Let us suppose now that $\alpha$ is a simple root which can be sent to $-\alpha_{\max}$
by means of an automorphism. Such a root has coefficient one in the decomposition of
$\alpha_{\max}$ and will be called a \emph{multiplicity free root}. Let us denote by
$\mathfrak{p}_{\alpha}$ the standard parabolic subalgebra of $\mathfrak{g}$ generated by
all root vectors corresponding to simple roots and their opposite except $-\alpha$. Let
$\mathfrak{l}_{\alpha}$ denote the set of all pairs in $\mathfrak{p}_{\alpha}\times
\mathfrak{p}_{\alpha}$ with equal Levi components. This is a Lagrangian subalgebra of
$\mathfrak{g}\times\mathfrak{g}$, where $\mathfrak{g}\times\mathfrak{g}$ has been endowed
with the following invariant bilinear form
\begin{equation}Q'((a,b),(c,d)):=K(a,c)-K(b,d).
\end{equation}

\begin{theorem}
The set of subalgebras $W\subseteq\mathbb{O}_{\alpha}\times\mathfrak{g}$, corresponding
to quasi-trigonometric solutions, is in a one-to-one correspondence with the set of
Lagrangian subalgebras $\mathfrak{l}$ of $\mathfrak{g}\times\mathfrak{g}$, with respect
to the bilinear form $Q'$, which satisfy the condition $\mathfrak{l}\oplus
\mathfrak{l}_{\alpha}=\mathfrak{g}\times\mathfrak{g}$.
\end{theorem}
\begin{proof}
The proof is based on the following result from \cite{S2}: Let $G$ be the simply
connected Lie group with Lie algebra $\mathfrak{g}$. Denote by $G_{\mathrm{ad}}$ the Lie group
$G/Z(G)$. Let $H$ be the Cartan subgroup with Lie algebra $\mathfrak{h}$ and $H_{\mathrm{ad}}$ its
image in $G_{\mathrm{ad}}$. If $\alpha$ is a multiplicity free root, then $\mathbb{O}_{\alpha}$
and $\mathbb{O}_{\alpha_{\max}}$ are conjugate by an element of
$H_{\mathrm ad}(\mathbb{C}((u^{-1})))$.

Suppose now that $W\subseteq\mathbb{O}_{\alpha}\times\mathfrak{g}$, then
$(\mathbb{O}_{\alpha}\times\mathfrak{g})^{\perp}\subseteq W^{\perp}=W$. It follows that
\begin{equation}
\frac{W}{(\mathbb{O}_{\alpha}\times\mathfrak{g})^{\perp}}\subseteq
\frac{\mathbb{O}_{\alpha}\times\mathfrak{g}}{(\mathbb{O}_{\alpha}\times\mathfrak{g})^{\perp}}\cong\frac{\mathbb{O}_{\alpha_{\max}}\times\mathfrak{g}}
{(\mathbb{O}_{\alpha_{\max}}\times\mathfrak{g})^{\perp}}\cong\mathfrak{g}\times\mathfrak{g}.
\end{equation}

Denote by $\mathfrak{l}$ the image of the quotient $\frac{W}{(\mathbb{O}_{\alpha}
\times\mathfrak{g})^{\perp}}$ in $\mathfrak{g}\times\mathfrak{g}$. One can check that
$\mathfrak{l}$ is a Lagrangian subalgebra of  $\mathfrak{g}\times\mathfrak{g}$ with
respect to $Q'$.

Moreover the image of $\mathfrak{g}[u]$ in $\mathfrak{g}\times\mathfrak{g}$, after
passing to the quotient as above, is precisely $\mathfrak{l}_{\alpha}$. Since $W$ is
transversal to $\mathfrak{g}[u]$, it follows that $\mathfrak{l}$ should be transversal to
$\mathfrak{l}_{\alpha}$.

Conversely, if $\mathfrak{l}$ is a Lagrangian subalgebra transversal to
$\mathfrak{l}_{\alpha}$ in  $\mathfrak{g}\times\mathfrak{g}$, then its corresponding
lift, $W$, is transversal to $\mathfrak{g}[u]$.
\end{proof}
We see that in the case of multiplicity free roots, the classification of
quasi-trigonometric solutions reduces to the following

\textbf{Problem.} Given a multiplicity free root $\alpha$, find all subalgebras
$\mathfrak{l}$ of $\mathfrak{g}\times\mathfrak{g}$ which build a Manin triple
$(Q',\mathfrak{l}_{\alpha},\mathfrak{l})$, with respect to the invariant bilinear form
$Q'$ on $\mathfrak{g}\times\mathfrak{g}$.

This problem was solved in \cite{P2} by using the classification of Manin triples for
complex reductive Lie algebras which had been obtained by P. Delorme in \cite{Del}. The
classification of Manin triples was expressed in terms of so-called \emph{generalized
Belavin-Drinfeld data}. Let us recall Delorme's construction.

Let $\mathfrak{r}$ be a finite-dimensional complex, reductive, Lie algebra and $B$ a
symmetric, invariant, nondegenerate bilinear form on $\mathfrak{r}$. The goal in
\cite{Del} is to classify all Manin triples of $\mathfrak{r}$ up to conjugacy under the
action on $\mathfrak{r}$ of the simply connected Lie group $\mathcal{R}$ whose Lie
algebra is $\mathfrak{r}$.

One denotes by $\mathfrak{r}_{+}$ and $\mathfrak{r}_{-}$ respectively the sum of the
simple ideals of $\mathfrak{r}$ for which the restriction of $B$ is equal to a positive
(negative) multiple of the Killing form. Then the derived ideal of $\mathfrak{r}$ is the
sum of $\mathfrak{r}_{+}$ and $\mathfrak{r}_{-}$.

Let $\mathfrak{j}_{0}$ be a Cartan subalgebra of $\mathfrak{r}$, $\mathfrak{b}_{0}$ a
Borel subalgebra containing $\mathfrak{j}_{0}$ and $\mathfrak{b}_{0}'$ be its opposite.
Choose $\mathfrak{b}_{0}\cap\mathfrak{r}_{+}$ as Borel subalgebra of $\mathfrak{r}_{+}$
and $\mathfrak{b}_{0}'\cap\mathfrak{r}_{-}$ as Borel subalgebra of $\mathfrak{r}_{-}$.
Denote by $\Sigma_{+}$ (resp., $\Sigma_{-}$) the set of simple roots of
$\mathfrak{r}_{+}$ (resp., $\mathfrak{r}_{-}$) with respect to the above Borel
subalgebras. Let $\Sigma=\Sigma_{+}\cup\Sigma_{-}$ and denote by $\mathcal{W}=
(H_{\alpha},X_{\alpha},Y_{\alpha})_{\alpha\in\Sigma}$ a Weyl system of generators of
$[\mathfrak{r},\mathfrak{r}]$. 
\begin{definition}[Delorme, \cite{Del}]\label{BD data}
One calls $(A,A',\mathfrak{i}_{\mathfrak{a}},\mathfrak{i}_{\mathfrak{a}'})$
\emph{generalized Belavin-Drinfeld data} with respect to $B$ when the following five
conditions are satisfied:

(1) $A$ is a bijection from a subset $\Gamma_{+}$ of $\Sigma_{+}$ on a subset
$\Gamma_{-}$ of $\Sigma_{-}$ such that
\begin{equation}\label{eq4}
B(H_{A\alpha},H_{A\beta})=-B(H_{\alpha},H_{\beta}), \alpha, \beta\in\Gamma_{+}.
\end{equation}

(2) $A'$ is a bijection from a subset $\Gamma'_{+}$ of $\Sigma_{+}$ on a subset
$\Gamma'_{-}$ of $\Sigma_{-}$ such that
\begin{equation}\label{eq4}
B(H_{A'\alpha},H_{A'\beta})=-B(H_{\alpha},H_{\beta}),\alpha,\beta\in\Gamma'_{+}.
\end{equation}

(3) If $C=A^{-1}A'$ is the map defined on $\rm{dom}(C)=\{\alpha\in\Gamma'_{+}:
A'\alpha\in\Gamma_{-}\}$ by $C\alpha=A^{-1}A'\alpha$, then $C$ satisfies:

For all $\alpha\in \rm{dom}(C)$, there exists a positive integer $n$ such that $\alpha$,...,
$C^{n-1}\alpha\in \rm{dom}(C)$ and $C^{n}\alpha\notin \rm{dom}(C)$.

(4) $\mathfrak{i}_{\mathfrak{a}}$ (resp., $\mathfrak{i}_{\mathfrak{a}'}$) is a complex
vector subspace of $\mathfrak{j}_{0}$, included and Lagrangian in the orthogonal
$\mathfrak{a}$ (resp., $\mathfrak{a}'$) to the subspace generated by $H_{\alpha}$,
$\alpha\in\Gamma_{+}\cup\Gamma_{-}$ (resp., $\Gamma'_{+}\cup\Gamma'_{-}$).

(5) If $\mathfrak{f}$ is the subspace of $\mathfrak{j}_{0}$ generated by the family
$H_{\alpha}+H_{A\alpha}$, $\alpha\in\Gamma_{+}$, and $\mathfrak{f}'$ is defined
similarly, then
\begin{equation}\label{eq5}
(\mathfrak{f}\oplus\mathfrak{i}_{\mathfrak{a}})\cap(\mathfrak{f}'\oplus
\mathfrak{i}_{\mathfrak{a}'})={0}.
\end{equation}
\end{definition}

Let $R_{+}$ be the set of roots of $\mathfrak{j}_{0}$ in $\mathfrak{r}$ which are linear
combinations of elements of $\Gamma_{+}$. One defines similarly $R_{-}$, $R'_{+}$ and
$R'_{-}$. The bijections $A$ and $A'$ can then be extended by linearity to bijections
from $R_{+}$ to $R_{-}$ (resp., $R'_{+}$ to $R'_{-}$). If $A$ satisfies condition (1),
then there exists a unique isomorphism $\tau$ between the subalgebra $\mathfrak{m}_{+}$
of $\mathfrak{r}$ spanned by $X_{\alpha}$, $H_{\alpha}$ and $Y_{\alpha}$,
$\alpha\in\Gamma_{+}$, and the subalgebra $\mathfrak{m}_{-}$ spanned by $X_{\alpha}$,
$H_{\alpha}$ and $Y_{\alpha}$, $\alpha\in\Gamma_{-}$, such that
$\tau(H_{\alpha})=H_{A\alpha}$, $\tau(X_{\alpha})=X_{A\alpha}$,
$\tau(Y_{\alpha})=Y_{A\alpha}$ for all $\alpha\in\Gamma_{+}$. If $A'$ satisfies (2), then
one defines similarly an isomorphism $\tau'$ between $\mathfrak{m'}_{+}$ and
$\mathfrak{m'}_{-}$.

\begin{theorem}[Delorme, \cite{Del}]\label{Manin}
(i) Let
$\mathcal{B}\mathcal{D}=(A,A',\mathfrak{i}_{\mathfrak{a}},\mathfrak{i}_{\mathfrak{a}'})$
be generalized Belavin-Drinfeld data, with respect to $B$. Let $\mathfrak{n}$ be the sum
of the root spaces relative to roots $\alpha$ of $\mathfrak{j}_{0}$ in
$\mathfrak{b}_{0}$, which are not in $R_{+}\cup R_{-}$. Let $\mathfrak{i}:=
\mathfrak{k}\oplus\mathfrak{i}_{\mathfrak{a}}\oplus\mathfrak{n}$, where
$\mathfrak{k}:=\{X+\tau(X):X\in\mathfrak{m}_{+}\}$.

Let $\mathfrak{n'}$ be the sum of the root spaces relative to roots $\alpha$ of
$\mathfrak{j}_{0}$ in $\mathfrak{b}_{0}'$, which are not in $R'_{+}\cup R'_{-}$. Let
$\mathfrak{i'}:=\mathfrak{k'}\oplus \mathfrak{i}_{\mathfrak{a'}}\oplus\mathfrak{n'}$,
where $\mathfrak{k'}:=\{X+\tau'(X):X\in\mathfrak{m'}_{+}\}$.

Then $(B,\mathfrak{i},\mathfrak{i}')$ is a Manin triple.

(ii) Every Manin triple is conjugate by an element of $\mathcal{R}$ to a unique Manin
triple of this type.
\end{theorem}
\begin{remark}
One says that the Manin triple constructed in (i) is associated to the generalized
Belavin-Drinfeld data $\mathcal{B}\mathcal{D}$ and the system of Weyl generators
$\mathcal{W}$. Such a Manin triple will be denoted by
$\mathcal{T}_{\mathcal{B}\mathcal{D},\mathcal{W}}$.
\end{remark}
Recall that $\Gamma$ denotes the set of simple roots relative to a Cartan subalgebra
$\mathfrak{h}$ of $\mathfrak{g}$. For a subset $S$ of $\Gamma$, let $[S]$ be the set of
roots in the linear span of $S$. Let
$\mathfrak{m}_{S}:=\mathfrak{h}+\sum_{\alpha\in[S]}\mathfrak{g}_{\alpha}$,
$\mathfrak{n}_{S}:=\sum_{\alpha>0,\alpha\notin [S]}\mathfrak{g}_{\alpha}$,
$\mathfrak{p}_{S}:=\mathfrak{m}_{S}+\mathfrak{n}_{S}$. We also consider
$\mathfrak{g}_{S}:=[\mathfrak{m}_{S},\mathfrak{m}_{S}]$,
$\mathfrak{h}_{S}:=\mathfrak{h}\cap\mathfrak{g}_{S}$ and
$\zeta_{S}:=\{x\in\mathfrak{h}:\alpha(x)=0, \forall\alpha\in S\}$.
Consider the Lagrangian subalgebra $\mathfrak{l}_{S}$ of $\mathfrak{g}\times\mathfrak{g}$
which consists of all pairs from $\mathfrak{p}_{S}\times\mathfrak{p}_{S}$ with equal
components in $\mathfrak{m}_{S}$. 

In \cite{P2}, the general result of Delorme was used in order to determine
Manin triples of the form $(Q',\mathfrak{l}_{S},\mathfrak{l})$. This enables one to
give the classification of all Lagrangian subalgebras $\mathfrak{l}$ of
$\mathfrak{g}\times\mathfrak{g}$ which are transversal to a given $\mathfrak{l}_{S}$.
We devote the rest of this section to presenting the main results of \cite{P2}. We refer to
\cite{P2} for the proofs. 

First of all, let us choose a suitable system of Weyl generators for
$\mathfrak{g}\times\mathfrak{g}$. Let $\mathfrak{b}$ be a Borel subalgebra of
$\mathfrak{g}$ containing the Cartan subalgebra $\mathfrak{h}$. Then
$\mathfrak{b}_{0}:=\mathfrak{b}\times\mathfrak{b}$ is a Borel subalgebra of
$\mathfrak{g}\times\mathfrak{g}$ containing the Cartan subalgebra
$\mathfrak{j}_{0}:=\mathfrak{h}\times\mathfrak{h}$. One denotes by $\Sigma_{+}$ the set
of pairs $(\alpha,0)$, and by $\Sigma_{-}$ the set of pairs $(0,-\alpha)$, where
$\alpha\in\Gamma$. Let $\Sigma:=\Sigma_{+}\cup\Sigma_{-}$.

Let $(X_{\alpha},Y_{\alpha},H_{\alpha})_{\alpha\in\Gamma}$ be an arbitrary Weyl system of
generators for $\mathfrak{g}$. A system of Weyl generators of
$\mathfrak{g}\times\mathfrak{g}$ with respect to this choice of simple roots can be
chosen as follows: $X_{(\alpha,0)}=(X_{\alpha},0)$, $H_{(\alpha,0)}=(H_{\alpha},0)$,
$Y_{(\alpha,0)}=(Y_{\alpha},0)$, $X_{(0,-\alpha)}=(0,Y_{\alpha})$,
$H_{(0,-\alpha)}=(0,-H_{\alpha})$, $Y_{(0,-\alpha)}=(0,X_{\alpha})$. A Manin triple
associated to some generalized Belavin-Drinfeld data $\mathcal{B}\mathcal{D}$ for
$\mathfrak{g}\times\mathfrak{g}$ with respect to this Weyl system will simply be noted by
$\mathcal{T}_{\mathcal{B}\mathcal{D}}$.

Let $\bar{\theta}_{S}$ be the automorphism of $\mathfrak{g}_{S}$ uniquely defined by the
properties $\bar{\theta}_{S}(X_{\alpha})=Y_{\alpha}$, $\bar{\theta}_{S}(Y_{\alpha})=
X_{\alpha}$ and $\bar{\theta}_{S}(H_{\alpha})=-H_{\alpha}$ for all $\alpha\in S$. Recall
that there is a short exact sequence
\begin{equation}
1\longrightarrow G_{S}\longrightarrow \mathrm{Aut}(\mathfrak{g}_{S})\longrightarrow \mathrm{Aut}_{S}
\longrightarrow 1
\end{equation}
where $\mathrm{Aut}_{S}$ denotes the group of automorphisms of the Dynkin diagram of
$\mathfrak{g}_{S}$. Let $\theta_{S}$ be the image of  $\bar{\theta}_{S}$ in $\rm{Aut}_{S}$.
Therefore $\bar{\theta}_{S}$ can be written uniquely as
\begin{equation}\label{theta}
\bar{\theta}_{S}=\psi_{S}\mathrm{Ad}_{g_{0}},
\end{equation}
where $g_{0}\in G_{S}$ and $\psi_{S}$ is the unique automorphism of $\mathfrak{g}_{S}$
satisfying the properties: $\psi_{S}(X_{\alpha})=X_{\theta_{S}(\alpha)}$,
$\psi_{S}(Y_{\alpha})=Y_{\theta_{S}(\alpha)}$,
$\psi_{S}(H_{\alpha})=H_{\theta_{S}(\alpha)}$ for all $\alpha\in S$.

\begin{theorem}\label{special Manin}
For any Manin triple $(Q',\mathfrak{l}_{S},\mathfrak{l})$, there exists a unique
generalized Belavin-Drinfeld data
$\mathcal{B}\mathcal{D}=(A,A',\mathfrak{i}_{\mathfrak{a}},\mathfrak{i}_{\mathfrak{a}'})$
where $A:S\times \{0\}\longrightarrow \{0\}\times(-S)$,
$A(\alpha,0)=(0,-\theta_{S}(\alpha))$ and $\mathfrak{i}_{\mathfrak{a}}=\rm{diag}(\zeta_{S})$,
such that $(Q',\mathfrak{l}_{S},\mathfrak{l})$ is conjugate to the Manin triple
$\mathcal{T}_{\mathcal{B}\mathcal{D}}=(Q',\mathfrak{i},\mathfrak{i'})$.

Moreover, up to a conjugation which preserves $\mathfrak{l}_{S}$, the Lagrangian
subalgebra $\mathfrak{l}$ is of the form \[\mathfrak{l}=(\mathrm{id}\times
\mathrm{Ad}_{g_{0}})(\mathfrak{i'}),\] where $g_{0}\in G_{S}$ is the unique element from the
decomposition (\ref{theta}).

\end{theorem}
\begin{lemma}\label{lemma1}
Let  $A:S\times \{0\}\longrightarrow \{0\}\times(-S)$,
$A(\alpha,0)=(0,-\theta_{S}(\alpha))$ and $\mathfrak{i}_{\mathfrak{a}}=\rm{diag}(\zeta_{S})$.
A quadruple $(A,A', \mathfrak{i}_{\mathfrak{a}},\mathfrak{i}_{\mathfrak{a'}})$ is
generalized Belavin-Drinfeld data if and only if the pair
$(A',\mathfrak{i}_{\mathfrak{a'}})$ satisfies the following conditions:

(1) $A':\Gamma_{1}\times \{0\} \longrightarrow \{0\}\times (-\Gamma_{2})$ is given by an
isometry $\tilde{A'}$ between two subsets $\Gamma_{1}$ and $\Gamma_{2}$ of $\Gamma$ :
$A'(\alpha,0)=(0,-\tilde{A'}(\alpha))$.

(2) Let $\rm{dom}(\tilde{A'},S):=\{\alpha\in\Gamma_{1}:\tilde{A'}(\alpha)\in S\cap
\Gamma_{2}\}$. Then for any $\alpha\in \rm{dom}(\tilde{A'},S)$ there exists a positive integer
$n$ such that $\alpha$, $(\theta_{S}\tilde{A'})(\alpha)$,...,
$(\theta_{S}\tilde{A'})^{n-1}(\alpha)\in \rm{dom}(\tilde{A'},S)$ but
$(\theta_{S}\tilde{A'})^{n}(\alpha)\notin \rm{dom}(\tilde{A'},S)$.

(3) Consider
$\Delta_{S}:=\{(h+h',-\psi_{S}(h)+h'):h\in\mathfrak{h}_{S},h'\in\zeta_{S}\}$.  Let
$\mathfrak{f'}$ be the subspace of $\mathfrak{h}\times\mathfrak{h}$ spanned by pairs
$(H_{\alpha},-H_{\tilde{A'}(\alpha)})$ for all $\alpha\in\Gamma_{1}$. Let
$\mathfrak{i}_{\mathfrak{a'}}$ be Lagrangian subspace of
$\mathfrak{a'}:=\{(h_{1},h_{2})\in\mathfrak{h}\times\mathfrak{h}:\alpha(h_{1})=0,
\beta(h_{2})=0, \forall\alpha\in\Gamma_{1}, \forall \beta\in\Gamma_{2}\}$. Then
\begin{equation}
(\mathfrak{f'}\oplus\mathfrak{i}_{\mathfrak{a'}})\cap\Delta_{S}={0}.
\end{equation}
\end{lemma}

\begin{definition}
A triple $(\Gamma_{1},\Gamma_{2},\tilde{A'})$ is called \emph{$S$-admissible} if
condition (2) of Lemma \ref{lemma1} is satisfied.
\end{definition}

Let $\Omega_{0}$ denote the Cartan component of the Casimir element $\Omega$. Let
$\pi_{1}$ (resp., $\pi_{2}$) be the projection of $\mathfrak{h}$ onto $\mathfrak{h}_{S}$
(resp., $\zeta_{S}$). Let $K_{0}$ be the restriction of the Killing form $K$ of
$\mathfrak{g}$ to $\mathfrak{h}$, which permits an identification between $\mathfrak{h}$
and $\mathfrak{h}^{*}$. If $R$ is an endomorphism of $\mathfrak{h}$, denote by $R^{*}$
the adjoint of $R$ regarded as an endomorphism of $\mathfrak{h}$.
\begin{lemma}\label{Lemma2}
(i) Suppose that $(\mathfrak{f'}\oplus\mathfrak{i}_{\mathfrak{a'}})\cap\Delta_{S}={0}$.
Then there exists a unique linear endomorphism $R$ of $\mathfrak{h}$ such that
\begin{equation}\label{eq7}
\mathfrak{f'}\oplus\mathfrak{i}_{\mathfrak{a'}}=\{(Rh,R'h):h\in\mathfrak{h}\},
\end{equation}
where $R'h:=\pi_{1}(h)-\pi_{2}(h)-\psi_{S}\pi_{1}(Rh)+\pi_{2}(Rh)$, and
\begin{equation}\label{eq8}
(\psi_{S}\pi_{1}+\pi_{2})R+R^{*}(\psi_{S}\pi_{1}+\pi_{2})=\mathrm{id}_{\mathfrak{h}}.
\end{equation}
(ii) There exists a bijection between the Lagrangian subspaces
$\mathfrak{i}_{\mathfrak{a'}}$ of $\mathfrak{a'}$ satisfying the condition
$(\mathfrak{f'}\oplus\mathfrak{i}_{\mathfrak{a'}})\cap\Delta_{S}={0}$, and the
endomorphisms $R$ of $\mathfrak{h}$ verifying (\ref{eq8}) and the additional condition:
\begin{equation}\label{eq9}
R((\psi_{S}\pi_{1}+\pi_{2})(H_{\gamma})+(\pi_{2}-\pi_{1})(H_{\tilde{A'}(\gamma)}))=
H_{\gamma},\forall\gamma\in\Gamma_{1}.
\end{equation}
(iii) There exists a bijection between endomorphisms $R$ of $\mathfrak{h}$ verifying
(\ref{eq8}) and (\ref{eq9}) and tensors $r\in\mathfrak{h}\otimes\mathfrak{h}$ satisfying
the following conditions:
\begin{equation}\label{eq10}
(\mathrm{id}\otimes(\psi_{S}\pi_{1}+\pi_{2}))(r)+((\psi_{S}\pi_{1}+\pi_{2})\otimes \mathrm{id})(r^{21})=
\Omega_{0},
\end{equation}
\begin{equation}\label{eq11}
((\tilde{A'}(\gamma)(\pi_{2}-\pi_{1})\otimes \mathrm{id})(r)=((\psi_{S}\pi_{1}+\pi_{2})\otimes
\gamma)(r),\forall\gamma\in\Gamma_{1}.
\end{equation}
\end{lemma}

\begin{corollary}
Let $A:S\times \{0\}\longrightarrow \{0\}\times(-S)$, $A(\alpha,0)=(0,-\theta_{S}
(\alpha))$ and $\mathfrak{i}_{\mathfrak{a}}=\rm{diag}(\zeta_{S})$. There exists a one-to-one
correspondence between generalized Belavin-Drinfeld data
$(A,A',\mathfrak{i}_{\mathfrak{a}},\mathfrak{i}_{\mathfrak{a'}})$ and pairs formed by an
$S$-admissible triple $(\Gamma_{1},\Gamma_{2},\tilde{A'})$ and a tensor
$r\in\mathfrak{h}\otimes\mathfrak{h}$ satisfying conditions (\ref{eq10}), (\ref{eq11}).
\end{corollary}

\begin{theorem}{\label{T}}
Suppose that $\mathfrak{l}$ is a Lagrangian subalgebra of
$\mathfrak{g}\times\mathfrak{g}$ transversal to $\mathfrak{l}_{S}$. Then, up to a
conjugation which preserves $\mathfrak{l}_{S}$, one has 
$\mathfrak{l}=(\mathrm{id}\times\mathrm{Ad}_{g_{0}})(\mathfrak{i'})$, where $\mathfrak{i'}$ is constructed from an $S$-admissible
triple $(\Gamma_{1},\Gamma_{2},\tilde{A'})$ and a tensor
$r\in\mathfrak{h}\otimes\mathfrak{h}$ satisfying conditions (\ref{eq10}), (\ref{eq11}).
\end{theorem}

Let $\alpha$ be a multiplicity free root of $\mathfrak{g}$. Set $S=\Gamma\setminus\{\alpha\}$. We write $\theta_{\alpha}$ instead of $\theta_{S}$ and $\psi_{\alpha}$ instead of $\psi_{S}$. We make the remark that a triple $(\Gamma_{1},\Gamma_{2},\tilde{A'})$ is $S$-admissible if and only if it is in one of the two situations:

I. If $\alpha\notin\Gamma_{2}$, then $(\Gamma_{1},\theta_{\alpha}(\Gamma_{2}),
\theta_{\alpha}\tilde{A'})$ is an admissible triple in the sense of \cite{BD1}.

II. If $\alpha\in\Gamma_{2}$ and $\tilde{A'}(\beta)=\alpha$, for some
$\beta\in\Gamma_{1}$, then $(\Gamma_{1}\setminus\{\beta\},
\theta_{\alpha}(\Gamma_{2}\setminus\{\alpha\}),\theta_{\alpha}\tilde{A'})$ is an
admissible triple in the sense of \cite{BD1}.

By applying Theorem \ref{T} in the particular case $S=\Gamma\setminus\{\alpha\}$ and
working with the tensor $\tilde{r}:=((\psi_{\alpha}\pi_{1}+\pi_{2})\otimes \mathrm{id})(r)$
instead of $r$, one obtains 
\begin{theorem}\label{maincor}
Let $\alpha$ be a multiplicity free root. Suppose that $\mathfrak{l}$ is a Lagrangian subalgebra of $\mathfrak{g}\times\mathfrak{g}$ transversal to $\mathfrak{l}_{\alpha}$. Then, up to a conjugation which preserves $\mathfrak{l}_{\alpha}$, one has $\mathfrak{l}=(\mathrm{id}\times \mathrm{Ad}_{g_{0}})(\mathfrak{i'})$, where $\mathfrak{i'}$ is constructed from a pair formed by $(\Gamma_{1},\Gamma_{2},\tilde{A'})$ and a tensor
$\tilde{r}\in\mathfrak{h}\otimes\mathfrak{h}$ satisfying the following conditions:

(1) $(\Gamma_{1},\Gamma_{2},\tilde{A'})$ is of type I or II from above.

(2) $\tilde{r}$ satisfies
\begin{equation}\label{eq12}
\tilde{r}+\tilde{r}^{21}=\Omega_{0},
\end{equation}

(3) If $(\Gamma_{1},\Gamma_{2},\tilde{A'})$ is of type I, then $\tilde{r}$ satisfies
\begin{equation}\label{eq13}
(\theta_{\alpha}\tilde{A'}(\gamma)\otimes \mathrm{id})(\tilde{r})+(\mathrm{id}\otimes\gamma)(\tilde{r})=0,\forall\gamma\in\Gamma_{1}.
\end{equation}

(4) If $(\Gamma_{1},\Gamma_{2},\tilde{A'})$ is of type II and $\tilde{A'}(\beta)=\alpha$,
then $\tilde{r}$ satisfies (\ref{eq13}) for all $\gamma\in\Gamma_{1}\setminus\{\beta\}$
and
\begin{equation}\label{eq14}
(\alpha(\pi_{2}-\psi_{\alpha}\pi_{1})\otimes \mathrm{id})(\tilde{r})=(\mathrm{id}\otimes\beta)(\tilde{r}).
\end{equation}
\end{theorem}
The construction of the quasi-trigonometric solutions can be summed up as follows. Suppose that
$(\Gamma_{1},\Gamma_{2},\tilde{A'})$ is of type I or II from above. Then one finds the
tensor $\tilde{r}$ and consequently $r$. This induces a unique endomorphism $R$ of
$\mathfrak{h}$, which in turn enables one to construct the subspace
$\mathfrak{i}_{\mathfrak{a'}}$, according to (\ref{eq7}). This is enough to reconstruct
$\mathfrak{i'}$. Then $\mathfrak{l}:=(\mathrm{id}\times \mathrm{Ad}_{g_{0}})(\mathfrak{i'})$ is a
Lagrangian subalgebra of $\mathfrak{g}\times\mathfrak{g}$ which is transversal to $\mathfrak{l}_{\alpha}$. Moreover, $\mathfrak{l}$ can
be lifted to a Lagrangian subalgebra $W$ of $\mathfrak{g}((u^{-1}))\oplus\mathfrak{g}$ which is transversal to $\mathfrak{g}[u]$.
By choosing two appropriate dual bases in $\mathfrak{g}[u]$ and  $W$ respectively, we
reconstruct the quasi-trigonometric solution $X(u,v)$. We will illustrate this procedure
by several examples.
\begin{example} \textbf{Quasi-trigonometric solutions for $\mathfrak{sl}(2)$.} Let ${e,f,h}$
be the canonical basis of $\mathfrak{sl}(2)$ and $\alpha$ be the simple root with root
vector $e$. Then $\Gamma=\{\alpha\}$. We have two cases:

I. $\Gamma_{1}=\Gamma_{2}=\emptyset$ and $\tilde{r}=\frac{1}{4} h\otimes h$.
Correspondingly we get one quasi-trigonometric solution:
\begin{equation}\label{trig1}
X_{0}(u,v)=\frac{v\Omega}{u-v}+r_{0},
\end{equation}
where $\Omega=e\otimes f+f\otimes e+\frac{1}{2}h\otimes h$ and $r_{0}=e\otimes f+
\frac{1}{4}h\otimes h$ is the Drinfeld-Jimbo r-matrix for $\mathfrak{sl}(2)$.

II.  $\Gamma_{1}=\Gamma_{2}=\{\alpha\}$, $\tilde{A'}=\rm id$ and
$\tilde{r}=\frac{1}{4}h\otimes h$. The corresponding quasi-trigonometric solution is
\begin{equation}\label{trig2}
X_{1}(u,v)=X_{0}(u,v)+(u-v)e\otimes e.
\end{equation}
\end{example}
\begin{example}
\textbf{Quasi-trigonometric solutions for $\mathfrak{sl}(3)$.} Denote by $\alpha$ the
simple root with root vector $e_{12}$ and by $\beta$ the one with root vector $e_{13}$.
Then $\Gamma=\{\alpha,\beta\}$ and both roots are singular. We will present the
quasi-trigonometric solutions corresponding to the root $\alpha$.

I. $\Gamma_{1}=\Gamma_{2}=\emptyset$. Then $\tilde{r}=a(e_{11}-e_{33})\otimes (e_{22}-
e_{33})+b(e_{22}-e_{33})\otimes (e_{11}-e_{33})+\frac{1}{3}(e_{11}-e_{33})\otimes
(e_{11}-e_{33})+\frac{1}{3}(e_{22}-e_{33})\otimes (e_{22}-e_{33})$, where
$a+b=-\frac{1}{3}$. The corresponding quasi-trigonometric solution is quasi-constant:
\begin{equation}\label{sol1}
X_{0}(u,v)=\frac{v\Omega}{u-v}+r_{0},
\end{equation}
where $r_{0}$ is the Drinfeld-Jimbo non-skewsymmetric r-matrix in $\mathfrak{sl}(3)$.

II. $\Gamma_{1}=\{\alpha\}$, $\Gamma_{2}=\{\alpha\}$, $\tilde{A'}(\alpha)=\alpha$. Then
$\tilde{r}=-\frac{1}{3}(e_{11}-e_{33})\otimes
(e_{22}-e_{33})+\frac{1}{3}(e_{11}-e_{33})\otimes
(e_{11}-e_{33})+\frac{1}{3}(e_{22}-e_{33})\otimes (e_{22}-e_{33})$. It follows that the
corresponding solution of the CYBE is again quasi-constant:
\begin{equation}\label{sol2}
X_{1}(u,v)=\frac{v\Omega}{u-v}+r_{1},
\end{equation}
where $r_{1}$ is another non-skewsymmetric r-matrix in $\mathfrak{sl}(3)$.

III.  $\Gamma_{1}=\{\alpha\}$, $\Gamma_{2}=\{\beta\}$, $\tilde{A'}(\alpha)=\beta$. Then
$\tilde{r}=-\frac{1}{3}(e_{22}-e_{33})\otimes (e_{11}-e_{33})
+\frac{1}{3}(e_{11}-e_{33})\otimes (e_{11}-e_{33})+\frac{1}{3}(e_{22}-e_{33})\otimes (e_{22}-e_{33})$. This data allows one to construct the following Lagrangian subalgebra which is transversal to $\mathfrak{l}_{\alpha}$: 
\begin{equation}
\mathfrak{l}=\{\left(\left(\begin{array}{ccc}
a&b&0\\
c&d&0\\
*&*&-a-d\end{array}\right),\left(\begin{array}{ccc}
-a-d&0&0\\
*&a&b\\
*&c&d\end{array}\right)\right): a,b,c,d\in\mathbb{C}\}.
\end{equation}

Correspondingly, one obtains the following solution:
\begin{equation}\label{sol3}
X_{2}(u,v)=X_{0}(u,v)-u(e_{12}\otimes e_{32})+v(e_{32}\otimes e_{12})-\frac{1}{6}(e_{11}-e_{22})\otimes (e_{22}-e_{33}).
\end{equation}

III'. $\Gamma_{1}=\{\beta\}$, $\Gamma_{2}=\{\alpha\}$, $\tilde{A'}(\beta)=\alpha$ and the
same $\tilde{r}$ as in III. We have a quasi-trigonometric solution which is gauge
equivalent to (\ref{sol3}).

IV.  $\Gamma_{1}=\Gamma_{2}=\{\alpha,\beta\}$, $\tilde{A'}(\alpha)=\beta$,
$\tilde{A'}(\beta)=\alpha$.  Then $\tilde{r}=-\frac{1}{3}(e_{22}-e_{33})\otimes
(e_{11}-e_{33}) +\frac{1}{3}(e_{11}-e_{33})\otimes
(e_{11}-e_{33})+\frac{1}{3}(e_{22}-e_{33})\otimes (e_{22}-e_{33})$. This data induces the following Lagrangian subalgebra: 
\[\mathfrak{l}=\{(X,TXT^{-1}):X\in\mathfrak{sl}(3)\},\]
where $T=e_{13}+e_{21}+e_{32}$. This is a subalgebra transversal to $\mathfrak{l}_{\alpha}$. 

The corresponding quasi-trigonometric solution is
\begin{equation}\label{sol4}
X_{3}(u,v)=X_{0}(u,v)-u(e_{12}\otimes e_{32}+e_{13}\otimes e_{12}+e_{12}\otimes e_{13})+ \end{equation}
\[v(e_{32}\otimes e_{12}+e_{12}\otimes e_{13}+e_{13}\otimes e_{12})+
(e_{13}+e_{23})\wedge e_{23}+\frac{1}{6}(e_{11}-e_{33})\wedge (e_{11}-e_{22}).\]

\begin{remark}
Solutions corresponding to the simple root $\beta$ are gauge equivalent to the solutions corresponding to $\alpha$. The solutions with non-constant polynomial part are the following: 
\begin{equation}
Y_{2}(u,v)=X_{0}(u,v)+v(e_{21}\otimes e_{23})-(ue_{23}\otimes e_{21})-\frac{1}{6}(e_{11}-e_{22})\otimes (e_{22}-e_{33}),
\end{equation}
which is equivalent to $X_{2}(u,v)$, and 
\begin{equation}
Y_{3}(u,v)=X_{0}(u,v)-u(e_{13}\otimes e_{23}+e_{23}\otimes e_{13}+e_{23}\otimes e_{21})+ \end{equation}
\[v(e_{13}\otimes e_{23}+e_{23}\otimes e_{13}+e_{21}\otimes e_{23})+
(e_{13}+e_{21})\wedge e_{12}+\frac{1}{6}(e_{11}-e_{22})\wedge (e_{22}-e_{33}),\]
which is equivalent to $X_{3}(u,v)$. 

\end{remark}
\end{example}

\begin{example}\textbf{Two examples of quasi-trigonometric solutions for $\mathfrak{sl}(N)$ of Cremmer-Gervais type.} We reconstruct two Lagrangian subalgebras which provide two quasi-trigonometric solutions
for $\mathfrak{g}=\mathfrak{sl}(N)$. These Lagrangian subalgebras are transversal to
$\mathfrak{l}_{\alpha_{1}}$ and are related to the Cremmer-Gervais Lie bialgebra
structure on $\mathfrak{g}$.

Let us consider the set of simple roots $\Gamma=\{\alpha_{1},...,\alpha_{N-1}\}$ and take
$S=\Gamma\setminus\{\alpha_{1}\}$. Let us denote by 
$(X_{\alpha},Y_{\alpha},H_{\alpha})_{\alpha\in\Gamma}$ the standard Weyl system.
In order to construct an $S$-admissible triple $(\Gamma_{1},\Gamma_{2},\tilde{A'})$, let
us first determine the map $\theta_{S}:S\longrightarrow S$. One can easily check that
$\theta_{S}$ is the following involution: $\theta_{S}(\alpha_{i})=\alpha_{N+1-i}$, for
all $i=2,...,N-1$, and that $g_{0}=T_{N-1}$ is the $(N-1)\times (N-1)$ matrix with 1 on the second diagonal and 0 elsewhere.

I. Consider $\Gamma_{1}=\{\alpha_{1},...,\alpha_{N-2}\}$, $\Gamma_{2}=\{\alpha_{2},...,
\alpha_{N-1}\}$ and $\tilde{A'}(\alpha_{i})=\alpha_{N-i}$. This is an $S$-admissible
triple. Indeed, $\theta_{S}\tilde{A'}(\alpha_{i})=\alpha_{i+1}$ and
$(\Gamma_{1},\theta_{S}(\Gamma_{2}),\theta_{S}\tilde{A'})$ is an admissible triple in the
sense of \cite{BD1}, which is known to be related to the Cremmer-Gervais Lie bialgebra
structure on $\mathfrak{g}$ (see \cite{CG}). The tensor $\tilde{r}$ satisfying
(\ref{eq12}), (\ref{eq13}) is the Cartan part of the Cremmer-Gervais non-skewsymmetric
constant r-matrix.

Let $\mathfrak{n}^{-}_{\alpha_{N-1}}$ denote the sum of all eigenspaces of negative roots
which contain $\alpha_{N-1}$ in their decomposition. Let $\mathfrak{n}^{-}_{\alpha_{1}}$ be
the sum of all eigenspaces of negative roots which contain $\alpha_{1}$ in their
decomposition. One can easily check that the Lagrangian subalgebra $\mathfrak{i'}$
constructed from this data is the following:
$\mathfrak{i'}=\mathfrak{n'}\oplus\mathfrak{i}_{\mathfrak{a'}}\oplus\mathfrak{k'}$, where
$\mathfrak{n'}=\mathfrak{n}^{-}_{\alpha_{N-1}}\times\mathfrak{n}^{-}_{\alpha_{1}}$,
$\mathfrak{i}_{\mathfrak{a'}}=span((diag(1,1,...,1,-N+1),(diag(-N+1,1,...,1))$ and $\mathfrak{k'}$ is spanned by the set of pairs $(X_{\alpha_{i}},Y_{\alpha_{N-i}})$,
$(Y_{\alpha_{i}},X_{\alpha_{N-i}})$, $(H_{\alpha_{i}},-H_{\alpha_{N-i}})$,
$i=1,...,N-2$. Let us consider $g_{0}$ as an element of $SL(N)$ and take
\begin{equation}
\mathfrak{l}_{1}=(\mathrm{id}\times \mathrm{Ad}_{g_{0}})(\mathfrak{i'})=\mathfrak{n}_{\alpha_{N-1}}^{-}\times\mathfrak{n}_{\alpha_{1}}^{-}\oplus\{(x,\tau(x)):x\in\mathfrak{m}_{\alpha_{N-1}}\}, 
\end{equation}
where $\tau(e_{ij})=e_{i+1,j+1}$ and $\mathfrak{m}_{\alpha_{N-1}}$ denotes the reductive part of $\mathfrak{p}_{\alpha_{N-1}}$. This Lagrangian subalgebra is
transversal to $\mathfrak{l}_{\alpha_{1}}$ and therefore induces a quasi-trigonometric solution corresponding to $\alpha_{1}$. Denote this solution by $X_{1}(u,v)$. 

II.  We consider $\Gamma_{1}=\Gamma_{2}=\Gamma$ and $\tilde{A'}(\alpha_{i})=\alpha_{N-i}$.
This is indeed an $S$-admissible triple since $(\Gamma_{1}\setminus\{\alpha_{N-1}\},
\theta_{S}(\Gamma_{2}\setminus\{\alpha_{1}\}),\theta_{S}\tilde{A'})$ is an admissible
triple in the sense of [1]. We have $\Gamma_{1}\setminus\{\alpha_{N-1}\}=
\{\alpha_{1},...,\alpha_{N-2}\}$,
$\theta_{S}(\Gamma_{2}\setminus\{\alpha_{1}\})=\{\alpha_{2},...,\alpha_{N-1}\}$ and
$\theta_{S}\tilde{A'}(\alpha_{i})=\alpha_{i+1}$. The tensor $\tilde{r}$, which satisfies
the system (\ref{eq12}), (\ref{eq13}) for $\gamma=\alpha_{i},i=1,...,N-2$ and
(\ref{eq14}), is as in case I. 

We obtain $\mathfrak{i}'=\{(X,\mathrm{Ad}_{U}(X)):X\in\mathfrak{sl}(N)\}$, 
where $U\in SL(N)$ is the matrix with 1 on the second diagonal and zero elsewhere. Finally take 
\begin{equation}
\mathfrak{l_{2}}=(\mathrm{id}\times \mathrm{Ad}_{g_{0}})(\mathfrak{i'})=\{(X,\mathrm{Ad}_{T}(X)):X\in\mathfrak{sl}(N)\},
\end{equation}
where $T=g_{0}U=e_{1N}+e_{21}+e_{32}+...+e_{N,N-1}$. This Lagrangian subalgebra is transversal
to $\mathfrak{l}_{\alpha_{1}}$ and consequently provides a quasi-trigonometric solution $X_{2}(u,v)$, corresponding again to the root $\alpha_{1}$.

Solutions $X_{1}(u,v)$ and $X_{2}(u,v)$ will be called \emph{quasi-trigonometric solutions of Cremmer-Gervais type}. Regarding their quantization, we will quantize instead the following solutions: 
\begin{equation}\label{Cremmer1}
X'_{1}(u,v)=(\sigma\otimes\sigma)X_{1}(u,v),
\end{equation}

\begin{equation}\label{Cremmer2}
X'_{2}(u,v)=(\sigma\otimes\sigma)X_{2}(u,v),
\end{equation}
where $\sigma(A)=-A^{t}$. 

\end{example}
\section{Quantum twists and their affinization }
The aim of the second part of our paper is to quantize certain quasi-trigonometric
solutions of the CYBE and the corresponding
Lie bialgebra stuctures on $\mathfrak{g} [u]$ in case $\mathfrak{g}=\mathfrak{sl}(N)$. We already know that all 
of them are in
the same twisting class and therefore the corresponding quantum groups are isomorphic
as algebras but with different comultiplications. However, these comultiplications
can be obtained from each other via quantum twisting (see \cite{Hal}, \cite{KPST}).

So, we  would like to outline some basic elements of quantum twisting  of
Hopf algebras (see \cite{ES}, p. 84-85). Suppose given a Hopf algebra $A:=A(m,\Delta,\epsilon,S$) with a
multiplication $m:A\otimes A\rightarrow A$, a coproduct $\Delta:A \rightarrow A\otimes
A$, a counit $\epsilon:A\rightarrow\mathbb C$, and an antipode $S:A\to A$. 

An invertible element $F\in A\otimes A$,
$F=\sum_i f^{(1)}_i\otimes f^{(2)}_i$ is called a quantum twist if it satisfies
the cocycle equation
\begin{equation}\label{qt2}
F^{12}(\Delta\otimes{\rm id})(F)=F^{23}({\rm id}\otimes\Delta)(F)\,,
\end{equation}
and the "unital" normalization condition
\begin{equation}\label{qt3}
(\epsilon \otimes{\rm id})(F)=({\rm id}\otimes\epsilon )(F)=1\,.
\end{equation}

Now we can define a twisted Hopf algebra
$A^{(F)}:=A^{(F)}(m,\Delta^{(F)},\epsilon,S^{(F)}$) which has the same multiplication $m$
and the counit mapping $\epsilon$ but the twisted coproduct and antipode
\begin{equation}\label{qt1}
\Delta^{(F)}(a)=F\Delta(a)F^{-1},\quad\;S^{(F)}(a)=u\,S(a)u^{-1}, \quad\;u= \sum_i
f^{(1)}_{i}S(f^{(2)}_i)\,\quad\; (a\in A)\,.
\end{equation} 

The Hopf algebra $A$ is called quasitriangular if it has an additional invertible element
(universal $R$-matrix) $R$, which relates the coproduct $\Delta$ with its
opposite coproduct $\tilde{\Delta}$ by the transformation
\begin{equation}\label{gt4}
\tilde{\Delta}(a)=R\,\Delta(a)R^{-1}\qquad (a\in A)~,
\end{equation}
with $R$ satisfying the quasitriangularity conditions
\begin{equation}\label{gt5}
(\Delta\otimes {\rm id})(R)=R^{13}R^{23}~,\qquad ({\rm id}\otimes\Delta)(R)\,=\,
R^{13}R^{12}~.
\end{equation}
The twisted ("quantized") Hopf algebra $A^{(F)}$ is also quasitriangular with the
universal $R$-matrix $R^{(F)}$ defined as follows
\begin{equation}\label{gt66}
R^{(F)} = F^{21}R\,F^{-1}~,
\end{equation}
where $F^{21}=\sum_i f^{(2)}_i\otimes f^{(1)}_{i}$.
So, the first step is to find a quantization of the Lie bialgebra structure on
$\mathfrak{g}[u]$ defined by $X_0(u,v)$, which was described in Section 3.
It is well-known that the corresponding quantum group is the so-called
$U_q(\mathfrak{g}[u])$ which is a parabolic subalgebra of the quantum affine algebra
$U_q (\hat{\mathfrak{g}})$. In case $\mathfrak{sl}(N)$ this algebra will be defined below.
However, it turns out that it is more convenient to work with its extended
version  $U_q(\mathfrak{gl}_N[u])$. 

The second step is to find explicit formulas for the quantum twists. We use two
methods.

{\it A). Affinization by Hopf isomorphism}. 
Let $F$ be a quantum twist and let $\mathop{Sup}(F)$ be a minimal
Hopf subalgebra, whose tensor square contains $F$, which we call the support
of $F$.
Similarly we define the support of a classical twist as the minimal Lie subbialgebra,
whose tensor square contains the given classical twist.

It turns out that for certain quasi-trigonometric solutions for $\mathfrak{sl}(N)$,
the corresponding support (in $\mathfrak{sl}(N)[u]!$) is isomorphic to the support
of a certain classical twist in $\mathfrak{sl}(N+1)$, which is however constant!
This observation enables us to apply results of \cite{ESS}, \cite{IO}, where constant twists
from the Belavin--Drinfeld list were quantized. Of course, the corresponding
quantum twists, one in  $U_q(\mathfrak{sl}(N)[u])$ and the second in $U_q(\mathfrak{sl}(N+1))$,
have isomorphic quantum supports. We will call this method {\it affinization by Hopf isomorphism.}

{\it B). Affinization by automorphism}. Let $F$ be some constant twist of
$U_q(\mathfrak{g})$ and $\omega$ be some automorphism of $U_q(\mathfrak{g}[u])$
such that $\omega (U_q(\mathfrak{g}))\nsubseteq U_q(\mathfrak{g})$. Then, under some conditions the element $F_{\pi(\omega)}:=(\omega^{-1}\pi(\omega)\otimes
\mathop{id})F$ will be also a quantum twist, i.e. it satisfies the cocycle equation
(\ref{qt2}). Here $\pi:U_q(\mathfrak{g}[u])\to U_q(\mathfrak{g})$  is the canonical
projection (the images of the affine roots are zero). The method is interesting on its own
but what is more important is that it leads to quantization of rational solutions
of the CYBE.

We consider these two methods on examples for the quantum algebra
$U_q(\mathfrak{gl}_N[u])$.

\section{A quantum seaweed algebra and its affine realization}

As we already mentioned it is more convenient to use instead of the simple Lie algebra
$\mathfrak{sl}_{N}^{}$ its central extension $\mathfrak{gl}_{N}^{}$. The polynomial
affine Lie algebra $\mathfrak{gl}_{N}[u]$ is generated by Cartan--Weyl basis
$e_{ij}^{(n)}:=e_{ij}u^n$ ($i,j=1,2,\ldots N$, $n=0,1,2,\ldots$) with the defining
relations
\begin{equation}\label{pr1}
[e_{ij}^{(n)},\,e_{kl}^{(m)}]=\delta_{jk}e_{il}^{(n+m)}-\delta_{il}e_{kj}^{(n+m)}~.
\end{equation}
The total root system $\Sigma$ of the Lie algebra $\mathfrak{gl}_{N}[u]$ with respect to
an extended Cartan subalgebra generated by the Cartan elements $e_{ii}$ ($i=1,2,\ldots,
N$) and $d=u(\partial/\partial u)$ is given by
\begin{equation}\label{pr2}
\mathop{\Sigma}(\mathfrak{gl}_{N}[u])\,=\,\{\epsilon_i-\epsilon_j,\,n\delta+
\epsilon_i-\epsilon_j,\,n\delta~|~i\neq j;~i,j=1,2,\ldots,N;~n=1,2,\ldots\}~,
\end{equation}
where $\epsilon_i^{}$ ($i=1,2,\ldots,N$) is the orthonormal basis of a $N$-dimensional
Euclidean space $\mathbb R^{N}$ dual to the Cartan subalgebra of $\mathfrak{gl}_N$.


We have
the following correspondence: $e_{ij}^{(n)}=e_{n\delta+\epsilon_i-\epsilon_j}$
for $i\neq j$, $n=0,1,2,\ldots$. 
We choose the following system of positive
simple roots:
\begin{equation}\label{pr3}
\mathop{\Pi}(\mathfrak{gl}_{N}[u])\,=\,\{\alpha_i:=\epsilon_i-\epsilon_{i+1},\,
\alpha_0:=\delta+\epsilon_N-\epsilon_1~|~i=1,2,\ldots N-1\}~.
\end{equation}
Now we would like to introduce seaweed algebra, which is important for our purposes.
Let $\mathfrak{sw}_{N+1}^{}$ be a subalgebra of $\mathfrak{gl}_{N+1}^{}$ generated by the
root vectors: $e_{21}^{}$, $e_{i,i+1}^{}$, $e_{i+1,i}^{}$ for $i=2,3,\ldots,N$ and
$e_{N,N+1}^{}$, and also by the Cartan elements: $e_{11}^{}+e_{N+1N+1}^{}$, $e_{ii}^{}$
for $i=2,3,\ldots,N$. It is easy to check that $\mathfrak{sw}_{N+1}$ has the structure of a
seaweed Lie algebra (see \cite{DK}).

Let $\hat{\mathfrak{sw}}_{N}^{}$ be a subalgebra of $\mathfrak{gl}_{N}^{}[u]$ generated
by the root vectors: $e_{21}^{(0)}$, $e_{i,i+1}^{(0)}$, $e_{i+1,i}^{(0)}$ for
$i=2,3,\ldots,N$ and $e_{N,1}^{(1)}$, and also by the Cartan elements: $e_{ii}^{(0)}$ for
$i=1,2,3,\ldots,N$. It is easy to check that the Lie algebras $\hat{\mathfrak{sw}}_{N}$
and $\mathfrak{sw}_{N+1}^{}$ are isomorphic. This isomorphism is described by the
following correspondence: $e_{i+1,i}^{}\leftrightarrow e_{i+1,i}^{(0)}$ for
$i=1,2,\ldots,N-1$, $e_{i,i+1}^{}\leftrightarrow e_{i,i+1}^{(0)}$ for $i=2,3,\ldots,N-1$,
$e_{N,N+1}^{}\leftrightarrow e_{N,1}^{(1)}$ for $i=2,3,\ldots,N-1$, and
$(e_{11}^{}+e_{N+1,N+1}^{})\leftrightarrow e_{11}^{(0)}$, $e_{ii}^{}\leftrightarrow
e_{ii}^{(0)}$ for $i=2,3,\ldots,N$. We shall call $\hat{\mathfrak{sw}}_{N}$ an \emph{affine
realization} of $\mathfrak{sw}_{N+1}^{}$.

Now let us consider the $q$-analogs of the previous Lie algebras. The quantum algebra
$U_{q}(\mathfrak{gl}_{N}^{})$ is generated by the Chevalley elements\footnote{We denote
the generators in the classical and quantum cases by the same letter "$e$". It should not
cause any misunderstanding.} $e_{i,i+1}^{}$, $e_{i+1,i}^{}$ $(i=1,2,\ldots,N-1)$, $q^{\pm
e_{ii}}$ $(i=1,2,\ldots,N)$ with the defining relations:
\begin{equation}\label{pr4}
\begin{array}{rcl}
q^{e_{ii}}q^{-e_{ii}}&=&q^{-e_{ii}}q^{e_{ii}}=1~,
\\[5pt]
q^{e_{ii}}q^{e_{jj}}&=&q^{e_{jj}}q^{e_{ii}}~,
\\[5pt]
q^{e_{ii}}e_{jk}^{}q^{-e_{ii}}&=&q^{\delta_{ij}-\delta_{ik}} e_{jk}^{}\quad(|j-k|=1)~,
\\[3pt]
[e_{i,i+1}^{},\,e_{j+1,j}^{}]&=&\delta_{ij}\,
\mbox{\large$\frac{q^{e_{ii}-e_{i+1,i+1}}\,-\,q^{e_{i+1,i+1}-e_{ii}}} {q\,-\,q^{-1}}$}~,
\\[5pt]
[e_{i,i+1}^{},\,e_{j,j+1}^{}]&=&0\quad{\rm for}\;\;|i-j|\geq 2~,
\\[5pt]
[e_{i+1,i}^{},\,e_{j+1,j}^{}]&=&0\quad{\rm for}\;\;|i-j|\geq 2~,
\\[5pt]
[[e_{i,i+1}^{},\,e_{j,j+1}^{}]_{q}^{},\,e_{j,j+1}^{}]_{q}^{}&=&0 \quad{\rm
for}\;\;|i-j|=1~,
\\[5pt]
[[e_{i+1,i}^{},\,e_{j+1,j}^{}]_{q}^{},\,e_{j+1,j}^{}]_{q}^{}&=&0 \quad{\rm
for}\;\;|i-j|=1~.
\end{array}
\end{equation}
where $[e_{\beta}^{},\,e_{\gamma}^{}]_{q}^{}$ denotes the $q$-commutator:
\begin{equation}\label{pr5}
[e_{\beta}^{},\,e_{\gamma}^{}]_{q}^{}\,:=\,e_{\beta}^{}e_{\gamma}^{}-
q^{(\beta,\gamma)}e_{\gamma}^{}e_{\beta}^{}~.
\end{equation}
The Hopf structure on $U_{q}(\mathfrak{gl}_{N}^{})$ is given by the following formulas
for a comultiplication $\Delta_{q}$, an antipode $S_{q}$, and a co-unit
$\varepsilon_{q}$:
\begin{equation}\label{pr6}
\begin{array}{rcl}
\Delta_{q}(q^{\pm e_{ii}})&=& q^{\pm e_{ii}}\otimes q^{\pm e_{ii}} ~,
\\[5pt]
\Delta_{q}(e_{i,i+1}^{}&=& e_{i,i+1}^{}\otimes 1+q^{e_{i+1,i+1}-e_{ii}}\otimes
e_{i,i+1}^{}~,
\\[5pt]
\Delta_{q}(e_{i+1,i}^{})&=& e_{i+1,i}^{}\otimes q^{e_{ii}-e_{i+1,i+1}}+1\otimes
e_{i+1,i}^{}~;
\end{array}
\end{equation}
\begin{equation}\label{pr7}
\begin{array}{rcl}
S_{q}(q^{\pm e_{ii}})&=&q^{\mp e_{ii}}~,
\\[5pt]
S_{q}(e_{i,i+1}^{})&=&-q^{e_{ii}-e_{i+1,i+1}}\,e_{i,i+1}^{}~,
\\[5pt]
S_{q}(e_{i+1,i}^{})&=&-e_{i+1,i}^{}\,q^{e_{i+1,i+1}-e_{i,i}}~; \phantom{aaaaaaaaaa}
\end{array}
\end{equation}
\begin{equation}\label{pr8}
\begin{array}{rcccl}
\phantom{aa}\varepsilon_{q}(q^{\pm e_{ii}})^{}&=&1~,\quad
\varepsilon_{q}(e_{ij}^{})&=&0\quad {\rm for}\;\;|i-j|=1~.
\end{array}
\end{equation}
In order to construct composite root vectors $e_{ij}^{}$ for $|i-j|\geq2$ we fix the
following normal ordering of the positive root system $\Delta_{+}^{}$ (see
\cite{T1,KT1,KT2})
\begin{equation}\label{pr9}
\begin{array}{c}
\epsilon_1^{}\!-\epsilon_2^{}\prec\epsilon_1^{}\!-\epsilon_3^{}
\prec\epsilon_2^{}\!-\epsilon_3^{}\prec\epsilon_{1}^{}\!-\epsilon_{4}^{}\prec
\epsilon_{2}^{}-\epsilon_{4}^{}\prec\epsilon_{3}^{}\!-\epsilon_{4}^{} \prec\ldots\prec
\\[5pt]
\epsilon_{1}^{}\!-\epsilon_{k}^{}\!\prec\epsilon_{2}^{}\!-\epsilon_{k}^{}
\!\prec\ldots\prec\epsilon_{k-1}^{}\!-\epsilon_{k}^{}\!\prec\ldots\prec
\epsilon_{1}^{}\!-\epsilon_{N}^{}\!\prec\epsilon_{2}^{}\!-\epsilon_{N}^{}
\!\prec\ldots\prec\epsilon_{N-1}^{}\!-\epsilon_{N}^{}~.
\end{array}
\end{equation}
According to this ordering we set
\begin{equation}\label{pr10}
e_{ij}^{}\,:=\,[e_{ik}^{},\,e_{kj}^{}]_{q^{-1}}^{},\qquad
e_{ji}^{}\,:=\,[e_{jk}^{},\,e_{ki}^{}]_{q}^{}~,
\end{equation}
where $1\le i<k<j\le N$. It should be stressed that the structure of the composite root
vectors does not dependent on the choice of the index $k$ in the r.h.s. of the definition
(\ref{pr10}). In particular, we have
\begin{equation}\label{pr11}
\begin{array}{rcccl}
e_{ij}^{}&:=&[e_{i,i+1}^{},\,e_{i+1,j}^{}]_{q^{-1}}^{}&= &[e_{i,j-1}^{},
\,e_{j-1,j}^{}]_{q^{-1}}^{}~,
\\[7pt]
e_{ji}^{}&:=&[\,e_{j,i+1}^{},\,e_{i+1,i}^{}\,]_{q}^{}&=&[e_{j,j-1}^{},\,
e_{j-1,i}^{}]_{q}^{}~,
\end{array}
\end{equation}
where $2\le i+1<j\le N$.

Using these explicit constructions and the defining relations (\ref{pr4}) for the
Chevalley basis it is not hard to calculate the following relations between
the Cartan--Weyl generators $e_{ij}$ ($i,j=1,2,\ldots, N$): 
\begin{eqnarray}
q^{e_{kk}^{}}e_{ij}^{}q^{-e_{kk}^{}}&=&q^{\delta_{ki}^{}-\delta_{kj}^{}}
e_{ij}^{}\qquad(1\le i,j,k\le N)~, \label{pr12}
\\[5pt]
[e_{ij}^{},\,e_{ji}^{}]&=&\frac{q^{e_{ii}^{}-e_{jj}^{}}-
q^{e_{jj}^{}-e_{ii}^{}}}{q-q^{-1}}\qquad(1\le i<j\le N)~, \label{pr13}
\\[5pt]
[e_{ij}^{},\,e_{kl}^{}]_{q^{-1}}&=&\delta_{jk}^{}e_{il}^{} \qquad(1\le i<j\le k<l\le N)~,
\label{pr14}
\\[5pt]
[e_{ik}^{},\,e_{jl}^{}]_{q^{-1}}^{}&=&(q-q^{-1})\,e_{jk}^{}e_{il}^{} \qquad(1\le
i<j<k<l\le N)~, \label{pr15}
\\[5pt]
[e_{jk}^{},\,e_{il}^{}]_{q^{-1}}^{}&=&0\qquad(1\le i\le j<k\le l\le N)~, \label{pr16}
\\[5pt]
[e_{kl}^{},\,e_{ji}^{}]&=&0\qquad(1\le i<j\le k<l\le N)~, \label{pr17}
\\[5pt]
[e_{il}^{},\,e_{kj}^{}]&=&0\qquad(1\le i<j<k<l\le N)~, \label{pr18}
\\[5pt]
[e_{ji}^{},\,e_{il}^{}]&=&e_{jl}^{}\,q^{e_{ii}^{}-e_{jj}^{}} \qquad(1\le i<j<l\le N)~,
\label{pr19}
\\[5pt]
[e_{kl}^{},\,e_{li}^{}]&=&e_{ki}^{}\,q^{e_{kk^{}}-e_{ll}^{}} \qquad(1\le i<k<l\le N)~,
\label{pr20}
\\[5pt]
[e_{jl}^{},\,e_{ki}^{}]&=&(q^{-1}-q)\,e_{kl}^{}e_{ji}^{}\, q^{e_{jj}-e_{kk}}\qquad(1\le
i<j<k<l\le N)~. \label{pr21}
\end{eqnarray}
These formulas can also be obtained from the relations between the elements of the Cartan--Weyl 
basis for  the
quantum superalgebra $U_q(\mathfrak{gl}(N|M)$ (see \cite{T3}). If we apply the Cartan
involution ($e_{ij}^{*}=e_{ji}^{}$) to the formulas above, we will get all relations between
elements of the Cartan--Weyl basis.

The quantum algebra $U_{q}(\mathfrak{gl}_{N}^{}[u])$ ($N\ge3$) is generated (as a unital
associative algebra) by the algebra
$U_{q}(\mathfrak{gl}_{N}^{})$ and the additional element $e_{N1}^{(1)}$ with the
relations:
\begin{equation}\label{pr22}
\begin{array}{rcl}
q^{\pm e_{ii}^{(0)}}e_{N1}^{(1)}&=&q^{\mp(\delta_{i1}-\delta_{iN})} e_{N1}^{(1)}q^{\pm
e_{ii}^{(0)}}~,
\\[5pt]
[e_{i,i+1}^{(0)},~e_{N1}^{(1)}]&=&0 \quad{\rm for}\;\, i=2,3,\ldots,N-2~,
\\[5pt]
[e_{i+1,i}^{(0)},~e_{N1}^{(1)}]&=&0\quad{\rm for}\;\, i=1,2,\ldots,N-1~,
\\[5pt]
[e_{12}^{(0)},~[e_{12}^{(0)},~e_{N1}^{(1)}]_{q}^{}]_{q}^{}&=&0~,
\\[5pt]
[e_{N-1,N}^{(0)},~[e_{N-1,N}^{(0)},~e_{N1}^{(1)}]_{q}^{}]_{q}^{}&=&0~,
\\[5pt]
[[e_{12}^{(0)},~e_{N1}^{(1)}]_{q}^{},~e_{N1}^{(1)}]_{q}^{}&=&0~,
\\[5pt]
[[e_{N-1,N}^{(0)},~e_{N1}^{(1)}]_{q}^{},~e_{N1}^{(1)}]_{q}^{}&=&0~.
\end{array}
\end{equation}
The Hopf structure of  $U_{q}(\mathfrak{gl}_{N}^{}[u])$ is defined by the formulas
(\ref{pr6})-(\ref{pr8}) for $U_{q}(\mathfrak{gl}_{N}^{(0)})$ and the following formulas
for the comultiplication and the antipode of $e_{N1}^{(1)}$:
\begin{eqnarray}\label{pr23}
\Delta_{q}(e_{N1}^{(1)})&=&e_{N1}^{(1)}\otimes 1+ q^{e_{11}^{(0)}-e_{NN}^{(0)}}\otimes
e_{N1}^{(1)}~,
\\[5pt]\label{pr24}
S_{q}(e_{N1}^{(1)})&=&-q^{e_{NN}^{(0)}-e_{11}^{(0)}}e_{N1}^{(1)}~.
\end{eqnarray}

Quantum analogs of the seaweed algebra $\mathfrak{sw}_{N+1}^{}$ and its affine
realization $\hat{\mathfrak{sw}}_{N}^{}$ are inherited from the quantum algebras
$U_{q}(\mathfrak{gl}_{N+1}^{})$ and $U_{q}(\mathfrak{gl}_{N}^{}[u])$. Namely, the quantum
algebra $U_q(\mathfrak{sw}_{N+1}^{})$ is generated by the root vectors: $e_{21}^{}$,
$e_{i,i+1}^{}$, $e_{i+1,i}^{}$ for $i=2,3,\ldots,N$ and $e_{N,N+1}^{}$, and also by the
$q$-Cartan elements: $q^{e_{11}^{}+e_{N+1,N+1}^{}}$, $q^{e_{ii}^{}}$ for $i=2,3,\ldots,N$
with the relations satisfying (\ref{pr4}). Similarly,  the quantum algebra
$U_q(\hat{\mathfrak{sw}}_{N}^{})$ is generated by the root vectors: $e_{21}^{(0)}$,
$e_{i,i+1}^{(0)}$, $e_{i+1,i}^{(0)}$ for $i=2,3,\ldots,N$ and $e_{N,1}^{(1)}$, and also
by the $q$-Cartan elements: $q^{e_{ii}^{(0)}}$ for $i=1,2,3,\ldots,N$ with the relations
satisfying (\ref{pr4}) and (\ref{pr22}). It is clear that the algebras
$U_{q}(\mathfrak{gl}_{N+1}^{})$ and $U_q(\hat{\mathfrak{sw}}_{N}^{})$ are isomorphic as
associative algebras but they are not isomorphic as Hopf algebras. However if we
introduce a new coproduct in the Hopf algebra $U_{q}(\mathfrak{gl}_{N+1}^{})$
\begin{equation}\label{pr25}
\Delta_{\;q}^{(\mathfrak{F}_{1,N+1}^{})}(x)\,=\,\mathfrak{F}_{1,N+1}^{}
\Delta_q(x)\mathfrak{F}_{1,N+1}^{-1}\quad (\forall x\in U_{q}(\mathfrak{gl}_{N+1}^{}))~,
\end{equation}
where
\begin{equation}\label{pr26}
\mathfrak{F}_{1,N+1}^{}:=q^{\,-e_{11}^{}\otimes e_{N+1,N+1}^{}}~,
\end{equation}
we obtain isomorphism of Hopf algebras
\begin{equation}\label{pr27}
U_q^{(\mathfrak{F}_{1,N+1})}(\mathfrak{sw}_{N+1}^{})\,\simeq\,
U_q(\hat{\mathfrak{sw}}_{N}^{})~.
\end{equation}
Here $U_q^{(\mathfrak{F}_{1,N+1})}(\mathfrak{sw}_{N+1}^{})$ denotes the
quantum seaweed algebra $U_q(\mathfrak{sw}_{N+1}^{})$ with the twisted coproduct
(\ref{pr25}).

\section{Cartan part of Cremmer-Gervais $r$-matrix}

First of all we recall classification of quasi-triangular $r$-matrices for a simple Lie
algebra $\mathfrak{g}$. The quasi-triangular $r$-matrices are solutions of the system
\begin{equation}\label{cg1}
\begin{array}{rcl}
r^{12}+r^{21}&=&\Omega~,
\\[7pt]
[r^{12},r^{13}]+[r^{12},r^{23}]+[r^{13},r^{23}]&=&0~,
\end{array}
\end{equation}
where $\Omega$ is the quadratic the Casimir two-tensor in
$\mathfrak{g}\otimes\mathfrak{g}$. Belavin and Drinfeld proved that any solution of this
system is defined by a triple $(\Gamma_{1},\Gamma_{2}, \tau)$, where
$\Gamma_{1},\Gamma_{2}$ are subdiagrams of the Dynkin diagram of $\mathfrak{g}$ and
$\tau$ is an isometry between these two subdiagrams. Further, each $\Gamma_i$ defines a
reductive subalgebra of $\mathfrak{g}$, and $\tau$ is extended to an isometry (with
respect to the corresponding restrictions of the Killing form) between the corresponding
reductive subalgebras of $\mathfrak{g}$. The following property of $\tau$ should be
satisfied: $\tau^{k}(\alpha)\not\in\Gamma_{1}$ for any $\alpha\in\Gamma_{1}$ and some
$k$. Let $\Omega_0$ be the Cartan part of $\Omega$. Then one can construct a
quasi-triangular $r$-matrix according to the following
\begin{theorem}[Belavin--Drinfeld \cite{BD1}]
Let $r_{0}\in\mathfrak{h}\otimes\mathfrak{h}$ satisfies the systems
\begin{eqnarray}\label{cg2}
r_{0}^{12}+r_{0}^{21}&=&\Omega_{0}~,
\\[5pt]\label{cg3}
(\alpha\otimes 1+1\otimes\alpha)(r_{0})&=&h_{\alpha}^{}
\\[5pt]\label{cg4}
(\tau(\alpha)\otimes 1+1\otimes\alpha)(r_{0})&=&0
\end{eqnarray}
for any $\alpha\in\Gamma_{1}$. Then the tensor
\begin{eqnarray}\label{cg5}
r=r_{0}+\sum_{\alpha>0}e_{-\alpha}^{}\otimes e_{\alpha}^{}+
\sum_{\alpha>0;k\geq1}e_{-\alpha}^{} \wedge e_{\tau^k(\alpha)}
\end{eqnarray}
satisfies (\ref{cg1}). Moreover, any solution of the system (\ref{cg1}) is of the above
form, for a suitable triangular decomposition of $\mathfrak{g}$ and suitable choice of a
basis $\{e_{\alpha}^{}\}$.
\end{theorem}
In what follows, aiming the quantization of algebra structures on the polynomial Lie
algebra $\mathfrak{gl}_{N}^{}[u])$ we shall use the twisted two-tensor $q^{r_{0}^{}(N)}$
where $r_{0}^{}(N)$ is the Cartan part of the Cremmer--Gervais $r$-matrix for the Lie
algebra $\mathfrak{gl}_{N}^{}$ when
$\Gamma_{1}=\{\alpha_{1},\alpha_{2},\ldots,\alpha_{N-2}\}$
$\Gamma_{2}=\{\alpha_{2},\alpha_{3},\ldots,\alpha_{N-1}\}$ and
$\tau(\alpha_i)=\alpha_{i+1}$. An explicit form of $r_{0}^{}(N)$ is defined by the
following proposition (see \cite{GG}).
\begin{proposition}
The Cartan part of the Cremmer--Gervais $r$-matrix for $\mathfrak{gl}_{N}^{}$ is given by
the following expression
\begin{equation}\label{cg6}
r_{0}^{}(\mathfrak{gl}_{N}^{})\,=\,
\displaystyle\frac{1}{2}~\sum_{i=1}^{N}e_{ii}^{}\otimes e_{ii}^{}+ \sum_{1\leq i<j\leq
N}\,\frac{N+2(i-j)}{2N}\;e_{ii}^{}\wedge e_{jj}^{}~.
\end{equation}
\end{proposition}
It is easy to check that the Cartan part (\ref{cg6}), $r_{0}^{}(N):=
r_{0}^{}(\mathfrak{gl}_{N}^{})$, satisfies the conditions
\begin{eqnarray}\label{cg7}
\bigl(\epsilon_k^{}\otimes{\rm id}+{\rm id}\otimes\epsilon_k^{}\bigr)
\bigl(r_{0}^{}(N)\bigr)&=&e_{kk}^{}\quad\;\;{\rm for}\;\;k=1,2,\ldots,N,
\\[5pt]\label{cg8}
\bigl(\epsilon_{k}^{}\otimes{\rm id}+{\rm id}\otimes\epsilon_{k'}^{}\bigr)
\bigl(r_{0}^{}(N)\bigr)&=&(k-k')\,\mathcal{C}_1^{}(N)-\!\!\!
\sum_{i=k'+1}^{k-1}e_{ii}^{}\quad\;{\rm for}\;\; 1\leq k'<k\leq N,
\end{eqnarray}
where $\mathcal{C}_1^{}(N)$ is the normalized central element (the Casimir element of
first order):
\begin{equation}\label{cg9}
\mathcal{C}_1^{}(N)\,:=\,\frac{1}{N}\sum_{i=1}^{N} e_{ii}^{}~.
\end{equation}
In particular (\ref{cg7}) and (\ref{cg8}) imply the Belavin--Drinfeld conditions
(\ref{cg3}) and (\ref{cg4}), i.e.
\begin{eqnarray}\label{cg10}
\bigl(\alpha_k^{}\otimes{\rm id}+{\rm id}\otimes\alpha_k^{}\bigr)
\bigl(r_{0}^{}(N)\bigr)&=& h_{\alpha_k}^{}\,:=\,e_{kk}^{}- e_{k+1,k+1}^{}~,
\\[5pt]\label{cg11}
\bigl(\tau(\alpha_{k'}^{})\otimes{\rm id}+{\rm id}\otimes\alpha_{k'}^{}\bigr)
\bigl(r_{0}^{}(N)\bigr)&=&\bigl(\alpha_{k'+1}^{}\otimes{\rm id}+ {\rm
id}\otimes\alpha_{k'}^{}\bigr)\bigl(r_{0}^{}(N)\bigr)\;=\;0
\end{eqnarray}
for $k=1,2\ldots,N-1~$, $k'=1,2\ldots,N-2~$, where $\alpha_k^{}=\epsilon_k^{}-
\epsilon_{k+1}^{}$ and $\alpha_{k'}^{}=\epsilon_{k'}^{}-\epsilon_{k'+1}^{}$ are the
simple roots of system $\mathop{\Pi}(\mathfrak{gl}_{N})$ (see (\ref{pr3})).

Now we consider some properties of the two-tensor $q^{r_{0}^{}(N)}$. First of all it is
obvious that this two-tensor satisfies cocycle equation. Further, for
construction of a quantum twist corresponding to the Cremmer--Gervais $r$-matrix
(\ref{cg6}) we introduce new Cartan--Weyl basis elements $e_{ij}^{\,\prime}$ ($i\neq j$)
for the quantum algebra $U_{q}(\mathfrak{gl}_{N})$ as follows
\begin{eqnarray}\label{cg12}
e_{ij}^{\,\prime}&=&\displaystyle e_{ij}^{}q^{\,\bigl((\epsilon_i- \epsilon_{j})\otimes
{\rm id}\bigr)\bigl(r_{0}(N)\bigr)}\,=
\,e_{ij}^{}\,q^{\;\sum\limits_{k=i}^{j-1}e_{kk}^{}-(j-i)\mathcal{C}_{1}^{}(N)}~,
\\[5pt]\label{cg13}
e_{ji}^{\,\prime}&=&q^{\,\bigl({\rm id}\otimes(\epsilon_{j}-
\epsilon_{i})\bigr)\bigl(r_{0}(N)\bigr)}\,e_{ji}^{}\,=\,
q^{\;\sum\limits_{k=i+1}^{j}e_{kk}^{}-(j-i)\mathcal{C}_{1}^{}(N)}\,e_{ji}^{}~,
\end{eqnarray}
for $1\leq i<j\leq N$. Permutation relations for these elements can be easily obtained from
the relations (\ref{pr12})--(\ref{pr21}). For example, we have
\begin{equation}\label{cg14}
\begin{array}{rcl}
[e_{ij}^{\prime},\,e_{ji}^{\prime}]&=&\displaystyle
[e_{ij}^{},\,e_{ji}^{}]\,q^{\,\bigl((\epsilon_i-\epsilon_{j})\otimes{\rm id}+{\rm
id}\otimes(\epsilon_{j}-\epsilon_{i})\bigr)\bigl(r_{0}^{}(N)\bigr)}
\\[12pt]
&=&\displaystyle\frac{q^{\,2\!\sum\limits_{k=i}^{j-1}e_{kk}^{}-
2(j-i)\mathcal{C}_{1}^{}(N)}-q^{\,2\!\!\!\sum\limits_{k=i+1}^{j}e_{kk}^{}-
2(j-i)\mathcal{C}_{1}^{}(N)}}{q-q^{-1}}~.
\end{array}
\end{equation}
It is not hard to check that the Chevalley elements $e_{i,i+1}^{\prime}$ and
$e_{i+1,i}^{\prime}$ have the following coproducts after twisting by the two-tensor
$q^{r_{0}^{}(N)}$:
\begin{equation}\label{cg15}
\begin{array}{rcl}
q^{\,r_{0}^{}(N)}\Delta_q^{}(e_{i,i+1}^{\prime})q^{-r_{0}^{}(N)}&=&
e_{i,i+1}^{\prime}\otimes q^{\,2\bigl((\epsilon_i-\epsilon_{i+1}) \otimes{\rm
id}\bigr)(r_{0}^{}(N))}+1\otimes e_{i,i+1}^{\prime}
\\[5pt]
&=&e_{i,i+1}^{\prime} \otimes q^{\,2e_{ii}^{}-2\mathcal{C}_{1}^{}(N)}+1\otimes
e_{i,i+1}^{\prime}~,
\end{array}
\end{equation}
\begin{equation}\label{cg16}
\begin{array}{rcl}
q^{\,r_{0}^{}(N)}\Delta_q^{}(e_{i+1,i}^{\prime})q^{-r_{0}^{}(N)}&=&
e_{i+1,i}^{\prime}\otimes1+q^{-2\bigl({\rm id}\otimes(\epsilon_{i+1}-
\epsilon_{i})\bigr)(r_{0}^{}(N))}\otimes e_{i+1,i}^{\prime}
\\[5pt]
&=&e_{i+1,i}^{\prime}\otimes1+ q^{2e_{i+1,i+1}^{}-2\mathcal{C}_{1}^{}(N)}\otimes
e_{i+1,i}^{\prime}~.
\end{array}
\end{equation}
for $1\leq i<N$. Since the quantum algebra $U_{q}(\mathfrak{gl}_{N})$ is a subalgebra of
the quantum affine algebra $U_{q}(\mathfrak{gl}_{N}[u])$ let us introduce the new affine
root vector $e_{N1}^{\prime(1)}$ in accordance with (\ref{cg12}):
\begin{eqnarray}\label{cg17}
e_{N1}^{\prime(1)}&=&\displaystyle e_{N1}^{}q^{\,\bigl((\epsilon_{N}^{}-
\epsilon_{1}^{})\otimes {\rm id}\bigr)\bigl(r_{0}(N)\bigr)}\,=
\,e_{N1}^{}\,q^{e_{NN}^{}-\mathcal{C}_{1}^{}(N)}~.
\end{eqnarray}
The coproduct of this element after twisting by the two-tensor $q^{r_{0}^{}(N)}$ has the
form
\begin{equation}\label{cg18}
\begin{array}{rcl}
q^{\,r_{0}^{}(N)}\Delta_q^{}(e_{N1}^{\prime(1)})q^{-r_{0}^{}(N)}&=&
e_{N1}^{\prime(1)}\otimes q^{\,2\bigl((\epsilon_{1}^{}-\epsilon_{N}^{}) \otimes{\rm
id}\bigr)(r_{0}^{}(N))}+1\otimes e_{N1}^{\prime(1)}
\\[7pt]
&=&e_{N1}^{\prime(1)}\otimes q^{\,2e_{NN}^{}-2\mathcal{C}_{1}^{}(N)}+ 1\otimes
e_{N1}^{\prime(1)}~.
\end{array}
\end{equation}

Consider the quantum seaweed algebra $U_{q}^{}(\mathfrak{sw}_{N+1}^{})$ after twisting by
the two-tensor $q^{r_{0}^{}(N+1)}$. Its new Cartan--Weyl basis and the coproduct for the
Chevalley generators are given by formulas (\ref{cg12}), (\ref{cg13}) and (\ref{cg15}),
(\ref{cg16}), where $N$ should be replaced  by $N+1$ and where $i\neq1$ in (\ref{cg12})
and (\ref{cg15}), and  $j\neq N$ in (\ref{cg13}), and $i\neq N$ in (\ref{cg16}). In
particular, for the element $e_{N,N+1}^{\prime}$ we have
\begin{eqnarray}\label{cg19}
e_{N,N+1}^{\prime}&=&e_{N,N+1}^{}\,q^{e_{NN}^{}- \mathcal{C}_{1}^{}(N+1)}~.
\\[5pt]
\label{cg20} q^{\,r_{0}^{}(N+1)}\Delta_q^{}(e_{N,N+1}^{\prime})q^{-r_{0}^{}(N+1)}&=
&e_{N,N+1}^{\prime}\otimes q^{\,2e_{NN}^{}-2\mathcal{C}_{1}^{}(N+1)}+ 1\otimes
e_{N,N+1}^{\prime}~.
\end{eqnarray}
Comparing the Hopf structures of the quantum seaweed algebra
$U_{q}^{}(\mathfrak{sw}_{N+1}^{})$  after twisting by the two-tensor
$q^{r_{0}^{}(N+1)}$and its affine realization $U_{q}^{}(\hat{\mathfrak{sw}}_{N}^{})$
after twisting by the two-tensor $q^{r_{0}^{}(N)}$ we see that these algebras are
isomorphic as Hopf algebras:
\begin{eqnarray}\label{cg21}
q^{\,r_{0}^{}(N+1)}\Delta_q^{}(U_{q}^{}(\mathfrak{sw}_{N+1}^{}))
q^{-r_{0}^{}(N+1)}&\simeq &q^{\,r_{0}^{}(N)}\Delta_q^{}
(U_{q}^{}(\hat{\mathfrak{sw}}_{N}^{}))q^{-r_{0}^{}(N)}~.
\end{eqnarray}
In terms of new Cartan--Weyl bases this isomorphism, $"\imath"$, is arranged as follows
\begin{eqnarray}\label{cg22}
\imath(e_{ij}^{\,\prime})&=&e_{ij}^{\prime(0)} \qquad{\rm for}\;\;2\leq i<j\leq N~,
\\[7pt]\label{cg23}
\imath(e_{ji}^{\,\prime})&=&e_{ji}^{\prime(0)} \qquad{\rm for}\;\;1\leq i<j\leq N-1~,
\\[7pt]\label{cg24}
\imath(e_{ii}^{}-\mathcal{C}_{1}^{}(N+1))&=&
e_{ii}^{(0)}-\mathcal{C}_{1}^{}(N)~,\qquad{\rm for}\;\;2\leq i\leq N~,
\\[2pt]\label{cg25}
\imath(e_{iN+1}^{\,\prime})\;=\;e_{i1}^{\prime(1)}&=&e_{i1}^{(1)}
q^{\,\bigl((\epsilon_{i}^{}-\epsilon_{1}^{})\otimes{\rm id}\bigr)
\bigl(r_{0}(N)\bigr)}\;=\;e_{i1}^{(1)}q^{\;\sum\limits_{k=i}^{N}
e_{kk}^{}-(N+1-i)\mathcal{C}_{1}^{}(N)}
\end{eqnarray}
for $2\leq i\leq N$, where the affine root vectors $e_{i1}^{(1)}$ ($2\leq i<N$) are
defined by the formula (cf. \ref{pr10}):
\begin{equation}\label{cg26}
e_{i1}^{(1)}\,=\,[e_{iN}^{(0)},\,e_{N1}^{(1)}]_{q^{-1}}^{}~.
\end{equation}

\section{Affine realization of Cremmer-Gervais twist}
In order to construct a twist corresponding to the Cremmer-Gervais $r$-matrix
(\ref{cg6}) we will follow the papers \cite{ESS},\cite{IO}.

Let $\mathcal{R}$ be a universal $R$-matrix of the quantum algebra
$U_{q}(\mathfrak{gl}_{N+1})$. According to \cite{KT1} it has the following form
\begin{equation}\label{cgt1}
\mathcal{R}\,=\,R\cdot K
\end{equation}
where the factor $K$ is a $q$-power of Cartan elements (see \cite{KT1}) and we do not
need its explicit form. The factor $R$ depends on the root vectors and it is given by the
following formula
\begin{equation}\label{cgt2}
\begin{array}{rcl}
R&=&R_{12}^{}(R_{13}^{}R_{23}^{})(R_{14}^{}R_{24}^{}R_{34}^{})
\cdots(R_{1,N+1}^{}R_{2,N+1}^{}\cdots R_{N,N+1}^{})
\\[5pt]
\!\!&=\!\!&\displaystyle\uparrow\prod_{j=2}^{N+1}\Bigl(\uparrow
\prod_{i=1}^{j-1}R_{ij}^{}\Bigr)~,
\end{array}
\end{equation}
where
\begin{equation}\label{cgt3}
R_{ij}^{}\,=\,\exp_{q^{-2}}((q-q^{-1})e_{ij}^{}\otimes e_{ji}^{})~,
\end{equation}
\begin{equation}\label{cgt4}
\exp_{q}(x):=\sum_{n\ge 0}\frac{x^{n}}{(n)_{q}!}~, \quad(n)_{q}!\equiv
(1)_{q}(2)_{q}\ldots (n)_{q}~,\quad (k)_{q}\equiv (1-q^{k})/(1-q)~.
\end{equation}
It should be noted that the product of factors $R_{ij}^{}$ in (\ref{cgt2}) corresponds to
the normal ordering (\ref{pr9}) where $N$ is replaced by $N+1$.

Let $R^{\,\prime}:=q^{r_{0}^{}(N+1)}R\,q^{-r_{0}^{}(N+1)}$. It is evident that
\begin{equation}\label{cgt5}
R^{\,\prime}\,=\,\uparrow\!\!\prod_{j=2}^{N+1}\Bigl(\uparrow\!
\prod_{i=1}^{j-1}R_{ij}^{\,\prime}\Bigr)~,
\end{equation}
where
\begin{equation}\label{cgt6}
R_{ij}^{\,\prime}\,=\,\exp_{q^{-2}}((q-q^{-1})e_{ij}^{\,\prime}\otimes
e_{ji}^{\,\prime})~.
\end{equation}
Here $e_{ij}^{\,\prime}$ and $e_{ji}^{\,\prime}$ are the root vectors (\ref{cg12}) and
(\ref{cg13}) where $N$ should be replaced by $N+1$.

Let $\mathcal{T}$ be a homomorphism which acts on the elements
$e_{ij}^{\,\prime}$ ($1\leq i<j\leq N+1$) by formulas
$\mathcal{T}(e_{ij}^{\,\prime})=e_{\tau(ij)}^{\,\prime}=e_{i+1,j+1}^{\,\prime}$ for
$1\leq i<j\leq N$, and $\mathcal{T}(e_{i,N+1}^{\,\prime})=0$ for all $i=1,2,\ldots, N$.
We set
\begin{equation}\label{cgt7}
R^{\,\prime(k)}\,:=\,(\mathcal{T}^k\otimes{\rm id})(R^{\,\prime})\,=\,\uparrow
\!\!\!\prod_{j=2}^{N+1-k}\Bigl(\uparrow\prod_{i=1}^{j-1}R_{ij}^{\,\prime(k)}\Bigr)~,
\end{equation}
where
\begin{equation}\label{cgt8}
R_{ij}^{\,\prime(k)}\,=\,\exp_{q^{-2}}\bigl((q-q^{-1})
\mathcal{T}^k(e_{ij}^{\,\prime})\otimes e_{ji}^{\,\prime}\bigr)\,=\,
\exp_{q^{-2}}\bigl((q-q^{-1})e_{i+k,j+k}^{\,\prime}\otimes e_{ji}^{\,\prime}\bigr)
\end{equation}
for $k\leq N-j$.

According to \cite{ESS}, \cite{IO}, the Cremmer-Gervais twist $\mathcal{F}_{CG}^{}$ in
$U_{q}(\mathfrak{gl}_{N+1})$ is given as follows
\begin{equation}\label{cgt9}
\mathcal{F}_{CG}^{}\,=\,F\cdot q^{r_{0}^{}(N+1)}~,
\end{equation}
where
\begin{equation}\label{cgt10}
F\,=\,R^{\,\prime(N-1)}R^{\,\prime(N-2)}\cdots R^{\,\prime(1)}~.
\end{equation}
It is easy to see that the support of the twist (\ref{cgt10}) is the
quantum seaweed algebra $U_{q}^{}(\mathfrak{sw}_{N+1}^{})$ with the coproducts
(\ref{cg15}) and (\ref{cg16}) where $N$ should be replaced by $N+1$. From the results of
the previous section it follows that we can immediately obtain an affine realization
$\hat{\mathcal{F}}_{CG}^{}$ which twists the quantum affine algebra
$U_{q}^{}(\mathfrak{gl}_{N}^{}[u])$:
\begin{eqnarray}\label{cgt11}
\hat{\mathcal{F}}_{CG}^{}&=&\hat{F}\cdot q^{r_{0}^{}(N)}~,
\\[7pt]\label{cgt12}
\hat{F}\,:=\,(\imath\otimes\imath)(F)&=&\hat{R}^{\,\prime(N-1)}
\hat{R}^{\,\prime(N-2)}\cdots\hat{R}^{\,\prime(1)}
\\[5pt]\label{cgt13}
\hat{R}^{\,\prime(k)}&=&\uparrow\!\!\!\prod_{j=2}^{N+1-k}
\Bigl(\uparrow\prod_{i=1}^{j-1}\hat{R}_{ij}^{\,\prime(k)}\Bigr)~.
\end{eqnarray}
where
\begin{eqnarray}\label{cgt14}
\hat{R}_{ij}^{\,\prime(k)}&=&\exp_{q^{-2}}\Bigl((q-q^{-1})
e_{i+k,j+k}^{\,\prime(0)}\otimes e_{ji}^{\,\prime(0)}\Bigr)\qquad {\rm for}\;\;1\leq
i<j\leq N-k~,
\\[5pt]\label{cgt15}
\hat{R}_{i,N+1-k}^{\,\prime(k)}&=&\exp_{q^{-2}}
\Bigl((q-q^{-1})e_{N,i+k}^{\,\prime(1)}\otimes e_{N+1-k,i}^{\,\prime(0)}\Bigr)\quad{\rm
for}\;\;1\leq i<N+1-k~.
\end{eqnarray}
Finally, using the isomorphism (\ref{pr27}) from Section 6 we obtain the following two results:
\begin{theorem}\label{CG1}
Let $\hat{\mathcal{F'}}_{CG}$ be the twist $\hat{\mathcal{F}}_{CG}$ reduced to $U_{q}(\hat{\mathfrak{sl}}_{N})$, 
and let $\hat{\cal{R}}$ be the universal $R$-matrix for $U_{q}(\hat{\mathfrak{sl}}_{N})$. 
Then the $R$-matrix $\hat{\mathcal{F'}}_{CG}^{21}\hat{\cal{R}}\hat{\mathcal{F'}}_{CG}^{-1}$ 
quantizes the quasi-trigonometric solution (\ref{Cremmer2}).
\end{theorem}

Now we turn to quantization of the quasi-trigonometric solution given by (\ref{Cremmer1}). 
The isomorphism (\ref{pr27}) shows that classical limits $\mathfrak{sw}_{N+1}$ and  $\hat{\mathfrak{sw}}_{N}$
are isomorphic as Lie bialgebras. Computations show that the support of solution (\ref{Cremmer1}) is
contained in the support of the solution (\ref{Cremmer2}). So, we can push the twist related to the solution 
(\ref{Cremmer1}) to $\mathfrak{sl}(N+1)$. It is not difficult to see that such an obtained twist will be defined by the
following Belavin-Drinfeld triple for  $\mathfrak{sl}(N+1)$: $\{\alpha_{2},...,\alpha_{N-1}\}\rightarrow\{\alpha_{3},...,\alpha_{N}\}$.
In fact, this is exactly the Cremmer-Gervais twist for $\mathfrak{sl}(N)$, embedded into $\mathfrak{sl}(N+1)$ as
the $N\times N$ block in the right low corner.

The corresponding constant twist for the above Belavin-Drinfeld triple can be quantized by means
of \cite{ESS}. Using once again (\ref{pr27}), we get a quantum twist which we denote by $\hat{\cal{F}}_{\cal{CG}}^{s}$.

\begin{theorem}\label{CG2}
Let $\hat{\cal{R}'}$ be the universal $R$-matrix for  $U_{q}(\hat{\mathfrak{sl}}_{N})$. 
Then the $R$-matrix $\hat{\cal{F}}_{\cal{CG}}^{s21}\hat{\cal{R}'}(\hat{\cal{F}}_{\cal{CG}}^{s})^{-1}$ 
quantizes the quasi-trigonometric solution (\ref{Cremmer1}).
\end{theorem}
\begin{remark}
In fact, using the isomorphism (\ref{pr27}) we can quantize all the quasi-trigonometric solutions of the 
CYBE corresponding to the first simple root $\alpha_1$ of $\mathfrak{sl}(N)$.

Let $W$ be the Lagrangian subalgebra of 
$\mathfrak{sl}(N)((u^{-1}))\oplus\mathfrak{sl}(N)$ 
contained in $\mathbb{O}_{\alpha_1}\oplus\mathfrak{sl}(N).$ Then it is not
difficult to show that the support of the corresponding
classical twist is contained in $\mathbb{O}_{\alpha_1}\cap\mathfrak{sl}(N)[u]$, which
is isomorphic to $\hat{\mathfrak{sw}}_{N}$. Therefore, it provides a 
classical twist in $\mathfrak{sw}_{N+1}$
since $\hat{\mathfrak{sw}}_{N}$ and $\mathfrak{sw}_{N+1}$ are isomorphic as Lie bialgebras.

Now we can again use results of \cite{ESS}, \cite{IO} to get the corresponding quantum affine twist.

\end{remark}

\section{Affinization by automorphism and quantization of rational $r$-matrices}

The aim of this section is to quantize certain rational $r$-matrices.
We begin with the following result:
\begin{theorem}
\label{affinizator} Let $\pi: U_{q}(\mathfrak{g}[u])\longrightarrow U_{q}(\mathfrak{g})$
be the canonical projection sending all the affine generators to zero.
Let $F\in
U_{q}(\mathfrak{g})\otimes U_{q}(\mathfrak{g})$ be a twist. 
Let us consider the following element
$$
F_{\pi\omega}=(\omega^{-1}\pi(\omega)\otimes 1)F(\pi\otimes{\rm
id})(\Delta(\omega^{-1}))\Delta(\omega)
$$
and represent it as a product $F'F$.
Then $F_{\pi\omega}$ is a twist iff
\begin{equation}
F_{12}( \pi\otimes{\rm id}\otimes{\rm id})(\Delta\otimes{\rm id})(F')\;=\;F'_{23}F_{12},
\label{rat2}
\end{equation} 
for some invertible $\omega\in U_{q}(\mathfrak{g}[u])$.
\end{theorem}
\begin{proof}
We will check the cocycle equation for an equivalent element
$$F_{\pi\omega}^{'}:=
(\omega\otimes\omega)F_{\pi\omega}\Delta(\omega^{-1})=(\pi\otimes{\rm
id})\left\{(\omega\otimes\omega)F\Delta(\omega^{-1})\right\}.$$ Note that
$F_{\pi\omega}^{'}\in U_{q}(\mathfrak{g})\otimes U_{q}(\mathfrak{g}[u])$. It follows that
$$
{\rm Assoc}(F_{\pi\omega}'):=(F_{\pi\omega}')_{12}(\Delta\otimes{\rm id})(F_{\pi\omega}')
({\rm id}\otimes\Delta)((F_{\pi\omega}')^{-1}) (F_{\pi\omega}')^{-1}_{23}\in
U_{q}(\mathfrak{g})\otimes U_{q}^{\otimes 2}(\mathfrak{g}[u]).
$$
On the other hand
$$
{\rm Assoc}(F_{\pi\omega}')= (\pi\otimes {\rm id}\otimes{\rm id})\left\{\omega^{\otimes
3}{\rm Assoc}(F_{\pi\omega})(\omega^{-1})^{\otimes 3}\right\}.
$$
If we take into account $ (\pi\otimes{\rm id})(F')=1\otimes 1$ and the property
(\ref{rat2}), then we get
\begin{equation}
\label{rat9}
\begin{array}{rcl}
&&{\rm Assoc}(F_{\pi\omega}')=
\\[5pt]
&&\quad={\rm Ad}(\pi(\omega)\otimes\omega^{\otimes 2})\bigr(F_{12}F'_{23}
(\Delta\otimes{\rm id})(F) ({\rm
id}\otimes\Delta)(F^{-1})(F_{23})^{-1}(F_{23}')^{-1}\bigr)
\\[7pt]
&&\quad={\rm Ad}(\pi(\omega)\otimes \omega^{\otimes 2})\bigl(F^{'}_{23}{\rm
Assoc}(F)(F^{'}_{23})^{-1}\bigr)~.
\end{array}
\end{equation}
Since $F$ is a twist  we deduce that ${\rm
Assoc}(F_{\pi\omega}')= 1\otimes 1\otimes 1$.\\
Conversely, let $F\in U_{q}(\mathfrak{g}[u])\otimes U_{q}(\mathfrak{g})$ and $(\pi\otimes
{\rm id})(F)$ is a twist, then there exist at least one $\omega\in
U_{q}(\mathfrak{g}[u])$ with the required property $(\ref{rat2})$. Indeed, note that
$(S\otimes S)(F_{21}^{-1})$ is a twist quantizing the same rational/quasi-trigonometric
$r-$matrix and thus there exists an invertible element $\omega$ such that
$$
(S\otimes S)(F_{21}^{-1})=(\omega\otimes \omega)F\Delta(\omega^{-1})\in
U_{q}(\mathfrak{g})\otimes U_{q}(\mathfrak{g}[u]).
$$
By taking projection $(\pi\otimes{\rm id})$ we obtain
$$
F=(\omega^{-1}\pi(\omega)\otimes 1)(\pi\otimes{\rm id})(F)(\pi\otimes{\rm
id})(\Delta(\omega^{-1}))\Delta(\omega).
$$
The property (\ref{rat2}) is necessary for the cocycle equation to hold.
\end{proof}
If an element $\omega$ satisfies the conditions of Theorem \ref{affinizator} we will call
it {\it{affinizator}} as it allows to construct an affine extension for a non-affine twist
$F$. Such element of course is not unique but some affinizators allow to construct
$F_{\pi\omega}$ which are compatible with the Yangian
degeneration.\\
Consider as an example the affinization of the coboundary twist
$$
F=(\exp_{q^{2}}(\lambda\mathop{}e_{-\alpha})\otimes \exp_{q^{2}}(\lambda\mathop{}
e_{-\alpha}))\Delta(\exp_{q^{-2}}(-\lambda\mathop{}e_{-\alpha}))
$$
with $\omega=\exp_{q^{2}}(\mu\mathop{}q^{-h_{\alpha}}e_{\delta-\alpha})$. In this case we
obtain
$$
F_{\pi\omega}=(1\otimes
1+(q^{2}-1)\lambda\mathop{}q^{-h_{\alpha}}e_{\delta-\alpha}\otimes
q^{-h_{\alpha}}+(q^{2}-1)\mu\mathop{}e_{-\alpha}\otimes 1)_{q^{2}}^{(-\frac 12 1\otimes
h_{\alpha})}.
$$
Let us form the equation for which the different $\omega$ are the solutions. In
order to find such $\omega$, we consider the following equation
$$
\mu({\rm id}\otimes
S)(F_{\pi\omega})=\sum_{i,j}\omega^{-1}\pi(\omega)F_{i}^{(1)}\pi(\omega^{(1)}_{j})S(\omega^{(2)}_{j})S(F_{i}^{2}),
$$
where $F=\sum_{i}F_{i}^{(1)}\otimes F_{i}^{(2)}$ and
$\Delta(\omega)=\sum_{i}\omega_{i}^{(1)}\otimes \omega_{i}^{(2)}$.

 Now we would like to explain how $\omega$-affinization can be used
to find a Yangian degeneration of the affine Cremmer--Gervais twists. Let us consider the
case $\mathfrak{sl}_3$. We set
\begin{equation}
\label{rat13} F={\cal F}_{CG_{3}}^{(\tau)}:=\exp_{q^{2}}(-(q-q^{-1})
\zeta\mathop{}\hat{e}_{12}^{(0)}\otimes \hat{e}_{32}^{(0)})\cdot\hat{{\cal K}_3}~,
\end{equation}
where
\begin{equation}\label{rat14}
\hat{\cal K}_3=q^{\frac 49 h_{12}\otimes h_{12}+\frac 29 h_{12}\otimes h_{23}+\frac 59
h_{23}\otimes h_{12}+ \frac 79 h_{23}\otimes h_{23}}
\end{equation}
with $h_{ij}:=e_{ii}-e_{jj}$. The twist (\ref{rat13}) belongs to
$U_{q}(\mathfrak{sl}_{3})\otimes U_{q}(\mathfrak{sl}_{3})[[\zeta]]$.

The following affinizator $\omega_{3}^{\rm long}$ was constructed in \cite{Sam1}. It is
given by the following formula
\begin{equation}\label{rat15}
\begin{array}{rcl}
\omega_{3}^{\rm long}\!\!&=\!\!&
\displaystyle\exp_{q^{2}}(\frac{\zeta}{1-q^{2}}\mathop{}q^{2h_{\alpha}^{\perp}}
\hat{e}_{21}^{(1)})\exp_{q^{2}}(-\frac{q\zeta^{2}}{(1-q^{2})^{2}}
q^{2h_{\beta}^{\perp}}\hat{e}_{31}^{(1)})
\\[9pt]
&&\times\,\displaystyle\exp_{q^{-2}}(\frac{\zeta^{2}}{1-q^{2}}
\mathop{}\hat{e}_{32}^{(0)})\exp_{q^{-2}}(\frac{\zeta}{1-q^{2}}\mathop{}
\hat{e}_{21}^{(0)})\exp_{q^{-2}}(\frac{\zeta^{2}}{1-q^{2}}\mathop{} \hat{e}_{32}^{(0)})~,
\end{array}
\end{equation}
where $h_{\alpha}^{\perp}=\frac 13 (e_{11}+e_{22})-\frac 23 e_{33}$ and
$h_{\alpha}^{\perp}=\frac 23 e_{11}- \frac 13 (e_{22}+e_{33})$.

For convenience sake we remind the reader that
\begin{equation}\label{rat16}
\begin{array}{lcl}
\hat{e}_{12}^{(0)}=e_{12}^{0}q^{h_{\beta}^{\perp}-h_{\alpha}^{\perp}}&&
\hat{e}_{21}^{(0)}=q^{h_{\beta}^{\perp}}e_{21}^{(0)}\\[2ex]
\hat{e}_{32}^{(0)}=q^{-h_{\beta}^{\perp}}e_{32}^{(0)},&&
\hat{e}_{31}^{0}=e_{32}^{0}e_{21}^{(0)}-q^{-1}e_{21}^{(0)}e_{32}^{(0)}\\[2ex]
\hat{e}_{31}^{(1)}=q^{h_{\alpha}^{\perp}-h_{\beta}^{\perp}}e_{31}^{(1)},&&
\hat{e}_{32}^{(1)}=e_{12}^{(0)}e_{31}^{(1)}-q e_{31}^{(1)}e_{12}^{(0)}\\[2ex]
\end{array}
\end{equation}
\begin{theorem}
The elements $\omega_{3}^{\rm long}$, $F={\cal F}_{CG_{3}}^{(\tau)}$ satisfy the
conditions of Theorem \ref{affinizator} and consequently $F_{\pi\omega}$ is a twist.
\end{theorem}
\begin{proof}
Straigtforward.
\end{proof}

It turns out that $F_{\pi\omega}'$ has a rational degeneration. To define this rational
degeneration we have to introduce the so-called $f$-generators:
\begin{equation}\label{rat17}
\begin{array}{lcl}
f_{0}=(q-q^{-1})\mathop{}\hat{e}_{31}^{(0)},&& f_{1}=q^{2h_{\beta}^{\perp}}
\hat{e}_{31}^{(1)}+q^{-1}\zeta\mathop{}\hat{e}_{31}^{(0)},\\[2ex]
f_{2}=(1-q^{-2})\mathop{}\hat{e}_{32}^{(0)},&& f_{3}=q^{h_{\alpha}^{\perp}}
\hat{e}_{32}^{(1)}-\zeta\mathop{}\hat{e}_{32}^{(0)}.
\end{array}
\end{equation}
Let us consider the Hopf subalgebra of  $U_{q}^{\hat{{\cal
K}}_3}(\hat{\mathfrak{sl}}_{3})$ generated by
$$\{h_{12}, h_{23}, f_{0}, f_{1},f_{2}, f_{3},
\hat{e}_{12}^{(0)},\hat{e}_{21}^{(0)}\}.$$ When $q\rightarrow 1$ we obtain the following
Yangian twist (see \cite{Sam1}):
\begin{equation}\label{rat18}
\begin{array}{rcl}
\overline{F}_{\pi\omega}' &=\!\!&(1\otimes1-\zeta\mathop{}1\otimes\overline{f}_{3}-
\zeta^{2}\mathop{}h_{\beta}^{\perp}\otimes \overline{f}_{2})^{(-h_{\beta}^{\perp}\otimes
1)}(1\otimes 1+ \zeta\mathop{}1\otimes\overline{e^{(0)}_{21}})
^{(-h_{\beta}^{\perp}\otimes 1)}
\\[7pt]
&&\times\;\exp(\zeta^{2}\mathop{}\overline{e^{(0)}_{12}}h_{13}
\otimes\overline{f}_{0})\exp(-\zeta\mathop{}
\overline{e^{(0)}_{12}}\otimes\overline{f}_{1})\cdot
\exp(-\zeta\mathop{}\overline{{e}_{12}^{(0)}}\otimes \overline{f}_{2})
\\[7pt]
&&\times\; (1\otimes 1-\zeta\mathop{}1\otimes\overline{f}_{3}-\zeta^{2}\mathop{}
h_{\alpha}^{\perp}\otimes\overline{f}_{2})^{((h_{\beta}^{\perp}-
h_{\alpha}^{\perp})\otimes1)}~,
\end{array}
\end{equation}
where the overlined generators are the generators of $Y(\mathfrak{sl}_{3})$. In the
evaluation representation we have:
\begin{equation}\label{rat19}
\begin{array}{ccc}
\overline{f}_{0}\mapsto e_{31}~,&& \overline{f}_{1}\mapsto u\mathop{}e_{31}
\\[5pt]
\overline{f}_{2}\mapsto e_{32}~, && {f}_{3}\mapsto u\mathop{}e_{32}
\\[5pt]
\overline{\vphantom{f}e}_{21}\mapsto e_{21}~,&& \overline{\vphantom{f}e}_{12}\mapsto
e_{12}~.
\end{array}
\end{equation}
Therefore we have obtained the following result:
\begin{theorem}
The Yangian twist $\overline{F}'_{\pi\omega}$ quantizes the following classical rational
$r-$matrix
\begin{equation}\label{rat20}
\begin{array}{rcl}
r(u,v)\!\!&=\!\!&\displaystyle\frac{\Omega}{u-v}+h_{\alpha}^{\perp}\otimes
ve_{32}-ue_{32} \otimes h_{\alpha}^{\perp}+h_{\beta}^{\perp}\wedge e_{21}
\\[10pt]
\!\!&&+e_{21}\otimes ve_{31}-ue_{31}\otimes e_{21}+e_{12}\wedge e_{32}~.
\end{array}
\end{equation}
\end{theorem}

To obtain a quantization of the second non-constant rational $r$-matrix for
$\mathfrak{sl}_{3}$ we take the following affinizator $\omega_{3}^{\rm\mathop{}short}$
and apply it to $F=q^{r_0(3)}$, where the Cartan part of the Cremmer-Gervais constant
$r$-matrix for $\mathfrak{sl}_{3}$ has the form:
\begin{equation}\label{rat21}
r_0(3)\;=\;\frac{2}{3}\bigr(h_{\alpha_1}\otimes h_{\alpha_1}+ h_{\alpha_2}\otimes
h_{\alpha_2}\bigl)+\frac{1}{3}\bigr(h_{\alpha_1}\otimes h_{\alpha_2}+h_{\alpha_2}\otimes
h_{\alpha_1}\bigl)+\frac{1}{6}h_{\alpha_1}\wedge h_{\alpha_2}~.
\end{equation}
We have
\begin{equation}\label{rat22}
\omega_{3}^{\rm short}\;=\;\exp_{q^{-2}}\bigl(\zeta\mathop{}\hat{e}_{21}^{(0)}
\bigr)\,\exp_{q^{2}}\Bigl(-\frac{\zeta}{1-q^{2}}\mathop{}q^{2h_{\alpha}^{\perp}}
\mathop{}\hat{e}_{31}^{(1)}\Bigr)\,\exp_{q^{-2}}\Bigl(\frac{\zeta}{1-q^{2}}\mathop{}
\hat{e}_{32}^{(0)}\Bigr)~,
\end{equation}
where
$$
\begin{array}{ccccc}
\hat{e}_{21}^{(0)}=q^{-\frac 13(h_{12}-h_{23})} e_{21}^{(0)},&&
\hat{e}_{32}^{(0)}=q^{-h_{\alpha}^{\perp}}e_{32}^{(0)},&&
\hat{e}_{31}^{(1)}=q^{-h_{\alpha}^{\perp}}e_{31}^{(1)}~.
\end{array}
$$
We have to calculate
\begin{equation}\label{rat23}
{\rm Aff}_{\omega^{\rm\mathop{} short}_3}(q^{r_0(3)})\;:= \;(\pi\otimes{\rm
id})\circ\Bigl((\omega_{3}^{\rm short}\otimes \omega_{3}^{\rm short})
q^{r_0(3)}\Delta(\omega_{3}^{\rm\ short})^{-1}\Bigr)~.
\end{equation}

Using standard commutation relations between $q$-exponents, $F'_{\pi\omega}$ can be
brought to the following form:
\begin{equation}\label{rat24}
\begin{array}{rcl}
&&\Bigr(1\otimes 1+\zeta\mathop{}1\otimes q^{2h_{\alpha}^{\perp}}
\hat{e}_{31}^{(1)}+\zeta\mathop{}q^{-2h_{\alpha}^{\perp}}\otimes\bigr({\rm
Ad}\exp_{q^{2}}(\zeta\hat{e}_{21}^{(0)})\bigr)(\hat{e}_{32}^{(0)})\Bigr)_{q^{2}}
^{(-h_{\alpha_{1}}^{\perp}\otimes 1)}
\\[7pt]
&&\qquad\times\;\Bigl(1\otimes 1+\zeta (1-q^{2})\mathop{}
1\otimes\hat{e}_{21}^{(0)}\Bigr)_{q^{-2}}^{(-\frac 13 (h_{12}-h_{23})\otimes
1)}q^{r_0(3)}~.
\end{array}
\end{equation}
The $q$-Hadamard formula allows us to calculate the Ad-term explicitly:
\begin{equation}\label{rat25}
\bigl({\rm Ad}\exp_{q^{-2}}(\zeta\mathop{}\hat{e}_{21}^{(0)})\bigr) (\hat{e}_{21}^{(0)})=
\hat{e}_{21}^{(0)}+\zeta\mathop{}q^{-h_{\beta}^{\perp}} e_{31}^{(0)}~,
\end{equation}
where $e_{31}^{(0)}:= e_{21}^{(0)}e_{32}^{(0)}-q\mathop{} e_{32}^{(0)}e_{21}^{(0)}$. To
define a rational degeneration we introduce $g$-generators, which satisfy the Yangian
relations as $q\to 1$:
$$
\begin{array}{ccccc}
g_{0}=(q-q^{-1})q^{-h_{\beta}^{\perp}}e^{(0)}_{31},&&
g_{1}=q^{2h_{\alpha}^{\perp}}\hat{e}_{31}^{(1)}+ \zeta\mathop{}q^{-h_{\beta}^{\perp}}
e^{(0)}_{31},&& g_{2}=(q^{2}-1)\mathop{}\hat{e}_{21}^{(0)}.
\end{array}
$$
Using $g$-generators we can calculate the rational degeneration of the twist
$F_{\pi\omega}'$:
\begin{equation}\label{rat26}
\begin{array}{rcl}
\overline{F}'_{\pi\omega}\!\!&=\!\!&\Bigl(1\otimes
1+\zeta\mathop{}1\otimes(\overline{g}_{1}+\overline{e^{(0)}_{32}})
-\zeta^2\mathop{}h_{\alpha}^{\perp}\otimes\overline{g}_{0}\Bigr)
^{(-h_{\alpha}^{\perp}\otimes 1)}
\\[7pt]&&\times\;\bigl(1\otimes 1-\zeta\mathop{}1\otimes
\overline{g}_{2}\bigr)^{(-\frac 13(h_{12}-h_{23})\otimes 1)}.
\end{array}
\end{equation}
\begin{theorem}
This Yangian twist $\overline{F}'_{\pi\omega}$ quantizes the following rational
$r$-matrix:
\begin{equation}\label{rat27}
r(u,v)\;=\;\frac{\Omega}{u-v}-u\mathop{}e_{31}\otimes
h_{\alpha}^{\perp}+v\mathop{}h_{\alpha}^{\perp}\otimes e_{31} +h_{\alpha}^{\perp}\wedge
e_{32}-\frac 13 (h_{12}-h_{23})\wedge e_{21}~.
\end{equation}
\end{theorem}
Therefore we have quantized all non-trivial rational $r$-matrices for $\mathfrak{sl}_3$
classified in \cite{S1}.
\section{Solutions for $\mathfrak{sl}(2)$ and deformed Hamiltonians}
We consider the case $\mathfrak{sl}(2)$. Let $\sigma^{+}=e_{12}$,
$\sigma^{-}=e_{21}$ and $\sigma^{z}=e_{11}-e_{22}$. 
Recall that in $\mathfrak{sl}(2)$ we have two quasi-trigonometric solutions, modulo gauge equivalence. The non-trivial solution is $X_{1}(z_{1},z_{2})=X_{0}(z_{1},z_{2})+(z_{1}-z_{2})(\sigma^{+}\otimes \sigma^{+})$. This solution is gauge equivalent to the following: 
\begin{equation}
X_{a,b}(z_{1},z_{2})=\frac{z_{2}\Omega}{z_{1}-z_{2}}+\sigma^{-}\otimes\sigma^{+}+\frac{1}{4}\sigma^{z}\otimes\sigma^{z}\label{E35}\end{equation}
\[
+a(z_{1}\sigma^{-}\otimes\sigma^{z}-z_{2}\sigma^{z}\otimes\sigma^{-})+b(\sigma^{-}\otimes\sigma^{z}-\sigma^{z}\otimes\sigma^{-}).\]

 The above quasi-trigonometric solution was quantized in \cite{KST2}. Let $\pi_{1/2}(z)$ be the two-dimensional vector
representation of $U_{q}(\widehat{sl_{2}}$). In this representation,
the generator $e_{-\alpha}$ acts as a matrix unit $e_{21}$, $e_{\delta-\alpha}$
as $ze_{21}$ and $h_{\alpha}$ as $e_{11}-e_{22}$. 
The quantum $R$-matrix of $U_{q}(\widehat{sl_{2}})$
in the tensor product $\pi_{1/2}(z_{1})\otimes\pi_{1/2}(z_{2})$ is
the following:\begin{equation}
R_{0}(z_{1},z_{2})=e_{11}\otimes e_{11}+e_{22}\otimes e_{22}+\frac{z_{1}-z_{2}}{q^{-1}z_{1}-qz_{2}}(e_{11}\otimes e_{22}+e_{22}\otimes e_{11})\label{E32}\end{equation}
\[
+\frac{q^{-1}-q}{q^{-1}z_{1}-qz_{2}}(z_{2}e_{12}\otimes e_{21}+z_{1}e_{21}\otimes e_{12}).\] 
\begin{proposition}
The $R$-matrix given by 
\begin{equation}
R:=R_{0}(z_{1},z_{2})+\frac{z_{1}-z_{2}}{q^{-1}z_{1}-qz_{2}}((b+az_{2})\sigma^{z}\otimes\sigma^{-}\label{E34}\end{equation}
\[
+(q^{-1}az_{1}+qb)\sigma^{-}\otimes\sigma^{z}+(b+az_{2})(q^{-1}az_{1}+qb)\sigma^{-}\otimes\sigma^{-})\]
is a quantization of the quasi-trigonometric solution $X_{a,b}$. 
\end{proposition}
\begin{corollary}
The rational degeneration  
\begin{equation}
R^{F}(u_{1},u_{2})=\frac{u_{1}-u_{2}}{u_{1}-u_{2}-\eta}(1-\eta\frac{P_{12}}{u_{1}-u_{2}}-\xi u_{2}\sigma^{z}\otimes\sigma^{-}\label{E41}\end{equation}
\[
+\xi(u_{1}-\eta)\sigma^{-}\otimes\sigma^{z}+\xi^{2}u_{2}(u_{1}-\eta)\sigma^{-}\otimes\sigma^{-}).\]
where $P_{12}$ denotes the permutation of factors in $\mathbb{C}^{2}\otimes\mathbb{C}^{2}$, 
is a quantization of the following rational solution of the CYBE: 
\begin{equation}
r(u_{1},u_{2})=\frac{\Omega}{u_{1}-u_{2}}+\xi(u_{1}\sigma^{-}\otimes\sigma^{z}-u_{2}\sigma^{z}\otimes\sigma^{-}).\label{E42}\end{equation}

\end{corollary}
The Hamiltonians of the periodic chains
related to the twisted $R$-matrix were computed in \cite{KST2}. We recall this result: 
We consider \begin{equation}
t(z)=Tr_{0}R_{0N}(z,z_{2})R_{0N-1}(z,z_{2})...R_{01}(z,z_{2})\label{ap1}\end{equation}
a family of commuting transfer matrices for the corresponding homogeneous
periodic chain, $[t(z'),$ $t(z")]=0$, where we treat $z_{2}$ as
a parameter of the theory and $z=z_{1}$ as a spectral parameter.
Then the Hamiltonian \begin{equation}
H_{a,b,z_{2}}=(q^{-1}-q)z\frac{d}{dz}t(z)\mid_{z=z_{2}}t^{-1}(z_{2})\label{ap2}\end{equation}
can be computed by a standard procedure: \begin{equation}
H_{a,b,z_{2}}=H_{XXZ}+\sum_{k}(C(\sigma_{k}^{z}\sigma_{k+1}^{-}+\sigma_{k}^{-}\sigma_{k+1}^{z})+D\sigma_{k}^{-}\sigma_{k+1}^{-}.\label{ap3}\end{equation}
Here $C=((q-1)/2)(b-az_{2}q^{-1})$, $D=(az_{2}+b)(q^{-1}az_{2}+qb)$,
$\sigma^{+}=e_{12}$, $\sigma^{-}=e_{21}$, $\sigma^{z}=e_{11}-e_{22}$
and \begin{equation}
H_{XXZ}=\sum_{k}(\sigma_{k}^{+}\sigma_{k+1}^{-}+\sigma_{k}^{-}\sigma_{k+1}^{+}+\frac{q+q^{-1}}{2}\sigma_{k}^{z}\sigma_{k+1}^{z}).\label{ap4}\end{equation}

We see that, by a suitable choice of parameters $a$, $b$ and $z_{2}$,
we can add to the XXZ Hamiltonian an arbitrary linear combination
of the terms $\sum_{k}\sigma_{k}^{z}\sigma_{k+1}^{-}+\sigma_{k}^{-}\sigma_{k+1}^{z}$
and $\sum_{k}\sigma_{k}^{-}\sigma_{k+1}^{-}$ and the model will remain
integrable.

Moreover, it was proved in \cite{KST2} that the Hamiltonian \begin{equation}
H_{\eta,\xi,u_{2}}=((q^{-1}-q)u-q^{-1}\eta)\frac{d}{du}t(u)\mid_{u=u_{2}}t^{-1}(u_{2})\label{ap5}\end{equation}
for \begin{equation}
t(u)=Tr_{0}R_{0N}(u,u_{2})R_{0N-1}(u,u_{2})...R_{01}(u,u_{2}),\label{ap6}\end{equation}
is given by the same formula (\ref{ap2}), where $C=\xi((q^{-1}-1)/2)u_{2}-(q^{-1}\xi\eta)/2$
and $D=\xi^{2}u_{2}(q^{-1}u_{2}-q\eta)$. Now it also makes sense
in the XXX limit $q=1$: \begin{equation}
H_{\eta,\xi,u_{2}}=H_{XXX}+\sum_{k}(C(\sigma_{k}^{z}\sigma_{k+1}^{-}+\sigma_{k}^{-}\sigma_{k+1}^{z})+D\sigma_{k}^{-}\sigma_{k+1}^{-}),\label{ap7}\end{equation}
where $C=-\xi\eta/2$ and $D=\xi^{2}u_{2}(u_{2}-\eta)$.

\section{Appendix}

In this appendix we give the proofs of the following results mentioned
in the text:

\begin{proposition}
Let $X$ be a rational or quasi-trigonometric solution of (\ref{eq:4}).
Then $X$ satisfies the unitarity condition (\ref{eq:5}). 
\end{proposition}
\begin{proof}
The proof is almost a word to word transcription of the proof of \cite{BD3},
Prop. 4.1. Interchanging $u_{1}$ and $u_{2}$ and also the first
and second factors in $\mathfrak{g}^{\otimes3}$ in equation (\ref{eq:4}),
we obtain \begin{equation}
[X^{21}(u_{2},u_{1}),X^{23}(u_{2},u_{3})]+[X^{21}(u_{2},u_{1}),X^{13}(u_{1},u_{3})]+\label{eq:B1}\end{equation}
\[
+[X^{23}(u_{2},u_{3}),X^{13}(u_{1},u_{3})]=0.\]
 Adding (\ref{eq:B1}) and (\ref{eq:4}), we get \begin{equation}
[X^{12}(u_{1},u_{2})+X^{21}(u_{2},u_{1}),X^{13}(u_{1},u_{3})+X^{23}(u_{2},u_{3})]=0.\label{eq:B2}\end{equation}

a) Suppose $X$ is rational, i.e. $X(u,v)=\frac{\Omega}{u-v}+p(u,v)$,
where $p$ is a polynomial. For $u_{1}$ and $u_{2}$ fixed, let us
multiply (\ref{eq:B2}) by $u_{2}-u_{3}$ and let $u_{3}\rightarrow u_{2}$.
It follows that \begin{equation}
[X^{12}(u_{1},u_{2})+X^{21}(u_{2},u_{1}),\Omega^{23}]=0.\label{eq:B3}\end{equation}
It is known that if a tensor $r\in\mathfrak{g}\otimes\mathfrak{g}$ satisfies
$[r\otimes1,\Omega^{23}]=0$, then $r=0$. It follows that $X^{12}(u_{1},u_{2})+X^{21}(u_{2},u_{1})=0$. 

b) Suppose $X$ is quasi-trigonometric, i.e. $X(u,v)=\frac{v\Omega}{u-v}+q(u,v)$
where $q$ is a polynomial function. By the same procedure we get
\begin{equation}
[X^{12}(u_{1},u_{2})+X^{21}(u_{2},u_{1}),u_{2}\Omega^{23}]=0\label{eq:B4}\end{equation}
which also implies the unitarity condition. 
\end{proof}
\begin{proposition}
Let $W$ be a Lie subalgebra satisfying conditions 2) and 3) of Theorem
\ref{thm:4}. Let $\widetilde{r}$ be constructed as in (\ref{eq:rtilda}).
Assume $\widetilde{r}$ induces a Lie bialgebra structure on $\mathfrak{g}[u]$
by $\delta_{\widetilde{r}}(a(u))=[\widetilde{r}(u,v),a(u)\otimes1+1\otimes a(v)]$.
Then $W\supseteq u^{-N}\mathfrak{g}[[u^{-1}]]$ for some positive $N$. 
\end{proposition}
\begin{proof}
Since $W$ is Lagrangian subalgebra, it is enough to prove that $W$
is bounded. Let us write

\begin{equation}
\widetilde{r}(u,v)=X_{0}(u,v)+\sum_{m}\Gamma_{m}\label{Ap1}\end{equation}
where $\Gamma_{m}$ is the homogeneous polynomial of degree $m$ with
coefficients in $\mathfrak{g}\otimes\mathfrak{g}$: \begin{equation}
\Gamma_{m}=\sum_{n+k=m}a_{mnk}u^{n}v^{k}.\label{Ap2}\end{equation}
It is enough to prove that there exists a positive integer $N$ such
that $\Gamma_{m}=0$ for $m\geq N$. 

We know that $\delta_{\widetilde{r}}(a)$ should belong to $\mathfrak{g}[u]\otimes\mathfrak{g}[v]$
for any element $a$ of $\mathfrak{g}$. On the other hand, one can
see that $[\Gamma_{m},a\otimes1+1\otimes a]$ is either 0 or has degree
$m$. This implies that $[\Gamma_{m},a\otimes1+1\otimes a]=0$ for
$m$ large enough. Therefore \begin{equation}
\Gamma_{m}=P_{m}(u,v)\Omega\label{Ap3}\end{equation}
 with $P_{m}(u,v)\in\mathbb{C}[[u,v]]$. Let us compute the following:

$[\Gamma_{m},au\otimes1+1\otimes av]=P_{m}(u,v)(u-v)[\Omega,a\otimes1]+P_{m}(u,v)v[\Omega,a\otimes1+1\otimes a]$\[
=P_{m}(u,v)(u-v)[\Omega,a\otimes1].\]

We choose an element $a$ such that $[\Omega,a\otimes1]\neq0$. We
obtain that if $P_{m}(u,v)$ is not identically zero then $P_{m}(u,v)(u-v)[\Omega,a\otimes1]$
is a homogeneous polynomial of degree $m+1$. Consequently, \begin{equation}
\delta_{\widetilde{r}}(au)=\sum_{m}P_{m}(u,v)(u-v)[\Omega,a\otimes1]\label{Ap4}\end{equation}
 cannot belong to $\mathfrak{g}[u]\otimes\mathfrak{g}[v]$ unless $P_{m}(u,v)=0$
for $m$ large enough. 
\end{proof}
\begin{theorem}
Let $X(z_1, z_2)$ be a quasi-trigonometric solution of the CYBE. 

1. Then there is a transformation
$\Psi(z)$, holomorphic around $z=1$, such that $$\left( \Psi(z_1)^{-1}\otimes \Psi(z_2)^{-1}\right)X(z_1,z_2)=Y(\frac{z_1}{z_2}),$$
where $Y(z)=\frac{\Omega}{z-1}+s(z)$, with $s(z)$ holomorphic around $z=1$.

2. $Y(e^{\lambda})$ is a trigonometric solution of the CYBE in the sense of Belavin--Drinfeld. 
\end{theorem}
\begin{proof}
Let us consider $X(z_1,z_2)=
\frac{z_2\Omega}{z_1-z_2}+p(z_1,z_2)$, where $p(z_1,z_2)$ is a
polynomial. Let $\{I_i\}$ be an orthonormal basis in $\mathfrak g$
with respect to the Killing form and $\{c_{ij}^k\}$ denote the
structure constants of $\mathfrak g$ with respect to $\{I_i\}$. Let
us write
$$p(z_1,z_2)= \sum_{i,j}p^{ij}(z_1,z_2)I_i\otimes I_j\ .$$
We set
$$h(z)= \sum_{i,j}p^{ij}(z,z)[I_i, I_j]\ 
=\sum_{i,j,k}p^{ij}(z,z)c_{ij}^kI_k\ .$$
Repeating the arguments of \cite{BD3}, one can prove that $h(z)$ and $X(z_1,z_2)$ satisfy 
$$z_1\frac{\partial X(z_1,z_2)}{\partial z_1}+
z_2\frac{\partial X(z_1,z_2)}{\partial z_2}=[h(z_1)\otimes 1+1\otimes 
h(z_2),
X(z_1,z_2)].$$
Suppose $\Psi(z)$ is a function with values in ${\mathrm {Aut}}({\mathfrak 
g})$,
which satisfies the differential equation
\begin{equation}
z\frac{d \Psi(z)}{d z}=(\mathrm{ad}h(z))\Psi(z),
\end{equation}
and the initial condition $\Psi(1)={\mathrm Id}$.
Then the function $Y(z_1,z_2)$ defined as
$$ Y(z_1,z_2)=\left( \Psi(z_1)^{-1}\otimes \Psi(z_2)^{-1}\right)X(z_1,z_2)\ 
,$$
satisfies the CYBE and depends on $z_1/z_2$ only. 

By construction,
$\Psi(z)$ is holomorphic in a neighborhood of $z=1$ and clearly 
$Y(z)=\frac{\Omega}{z-1}+s(z)$, where $s(z)$ is a holomorphic function in the same neighborhood.

Now we turn to the proof of the second statement.
Apply the change of variables $z_1=e^u$, $z_2=e^v$. Then let us prove that $\Psi(e^u)$ is 
holomorphic
in the entire complex plane and $Y(u,v)$ is a trigonometric solution of the 
CYBE.
 
Let $\Psi_1(u):=\Psi(e^u)$. Clearly this operator satisfies the equation 
\begin{equation}\label{diff eq}
\frac{d \Psi_1(u)}{du}=(\mathrm{ad}h_1(u))\Psi_1(u),
\end{equation}
where $h_1(u)=h(e^u)$. 

Let $U(u)$  be the matrix of the operator $\Psi_1(u)$ in the basis $\{ I_{i}\}$.
Let $a^{ij}(u,v):=p^{ij}(e^u,e^v)$ (the decomposition of $h_1(u)$ in the basis $\{ I_{i}\}$). 
Equation (\ref{diff eq}) is equivalent to 
\begin{equation}
\frac{dU(u)}{du}=H(u)U(u),\label{eq:matrix}
\end{equation}
where $H(u)$ is the matrix with elements 
\begin{equation}
h_{ij}(u)=\sum_{s,r,t}c_{sj}^{i}c_{rt}^{s}a^{rt}(u,u).
\label{Ap13}\end{equation}

Since the matrix function $H(u)$ is holomorphic in $\mathbb{C}$, the matrix
equation (\ref{eq:matrix}) admits a unique solution satisfying $U(0)=E$.
This solution is holomorphic in $\mathbb{C}$ 
because $U(u)=P\rm{exp}(\int_{0}^{u}H(v)dv)$
(ordered exponential) and \[
\left\Vert U(u)\right\Vert =\left\Vert 1+\int_{0}^{u}H(v)dv+\int_{0}^{u}(\int_{0}^{v_{1}}H(v_{1})H(v_{2})dv_{2})dv_{1}+...\right\Vert \leq\]
\[
\leq\exp(\int_{0}^{u}\left\Vert H(v)\right\Vert dv.\]

Moreover, according to \cite{BD3}, the linear operator $\Psi(u)$, corresponding
to $U(u)$, is an automorphism of $\mathfrak{g}$. 

Clearly $Y(e^u,e^v)$ depends only on $u-v$ and, as a function in one variable, has poles when $e^{u-v}=1$. Hence 
it is trigonometric in the sense of Belavin--Drinfeld. 
This ends the proof.
\end{proof} 
\subsection*{Acknowledgments}
The paper has been partially supported by the Royal Swedish Academy of Sciences under the
program "Cooperation between Sweden and former USSR" and the grants RFBR-05-01-01086,
INTAS-OPEN-03-51-3350 (V.N.T.).

The authors are thankful to the referee of the first version of the paper for
pointing out the importance of the paper \cite{Del}.

\end{document}